\documentstyle[11pt]{article}
\def\bold{\bf}
\def\Bbbd{\bf}
\def\be{\begin{equation}}
\def\ee{\end{equation}}
\def\bear{\begin{eqnarray}}
\def\eear{\end{eqnarray}}
\def\best{\begin{eqnarray*}}
\def\eest{\end{eqnarray*}}
\def\pf{{\bf Proof}: }
\renewcommand{\theequation}{\arabic{section}.\arabic{equation}}
\newtheorem{theorem}{Theorem}[section]
\newtheorem{prop}[theorem]{Proposition}

\newtheorem{lemma}[theorem]{Lemma}
\newtheorem{cor}[theorem]{Corollary}
\newtheorem{defn}[theorem]{Definition}

\newtheorem{fact}[theorem]{Fact}
\def\rem{ \addtocounter{theorem}{1}
{\non \bf Remark \arabic{section}.\arabic{theorem} }}
\def\ker{{\rm Ker D}}
\def\cok{{\rm Coker D}}

\def\step#1{\vskip.1in  \non{\bf #1.} }  
\def\non{\noindent}
\def\pf{\non {\bf Proof. }}
\def\qed{\nopagebreak \hskip .1in { $\Box$ }\penalty10000 \hskip\parfillskip
\par  }
\def\ra{\rightarrow}
\def\rg{\rangle}
\def\rgl{\rg_{\vphantom{I}_\la}}
\def\lg{\langle}
\def\r){\right)}
\def\l({\left(}
\def\ma(#1){\mathop {#1} \limits}
\def\fl{{f_\la}}

\def\b{\beta}
\def\bl{\b_\la}
\def\Si{\Sigma}
\def\De{\Delta}
\def\ti{\times}

\def\del{\overline \partial_{J}}
\def\mspa{{\vphantom{\ma(\int)_{k}}}}
\def\Z{{ \Bbbd Z}}
\def\R{{ \Bbbd R}}
\def\P{{ \Bbbd P}}

\def\cx{{ \Bbbd C}}
\def\pr{\pi_{-}}

\def\Dl{D_\la}
\def\Dls{D^*_\la}
\def\La{\Lambda^{0,1}}
\def\Lae{\Lambda^{0,1}_E}
\def\Lal{\Lambda^{0,1}_{low}}
\def\Lo{\Lambda^{0}}
\def\Lol{\Lambda^{0}_{low}}
\def\Loe{\Lambda^{0}_{E}}
\def\ho(#1){ H^{1}( T^2, {#1}^{*}T \P^n)}
\def\lo(#1){\Lambda^{0}({#1}^{*}T \P^n)}
\def\lol(#1){\Lambda^{0}_{low}({#1}^{*}T \P^n)}
\def\loe(#1){\Lambda^{0}_{E}({#1}^{*}T\P^n)}
\def\lla(#1){\Lambda^{0,1}({#1}^{*}T\P^n)}
\def\lal(#1){\Lambda^{0,1}_{low}({#1}^{*}T\P^n)}
\def\lae(#1){\Lambda^{0,1}_E ({#1}^{*}T\P^n)}

\def\lp{{\|}_{p,\la}}
\def\lpi{{\|}_{1,p,\la}}
\def\lu{{\|}_{2,\la}}
\def\lui{{\|}_{1,2,\la}}
  
\def\la{\lambda}
\def\sla{\sqrt \lambda}
\def\al{\alpha}
\def\o{ O} 
\def\co{{\cal O}}
\def\n{\nu}
\def\ze{\zeta}

\def\tel{\theta_{\la}}

\def\de{\delta}
\def\ep{\varepsilon}
\def\bz{{z}}
\def\bn{{n}}
\def\cmo{{\cal M}^0}

\def\cm{{\cal M}}

\def\cn{{\cal N}}
\def\cz{{\cal Z}}
\def\czh{\widehat{\cal Z}}
\def\cw{{\cal W}}
\def\czt{\widetilde{\cal Z}}

\def\cu{{\cal U}}

\def\cut{\widetilde{\cal U}}
\def\cud{{\bf U}_\delta}
\def\oml{\omega_{\la}}
\def\om{\omega}
\def\omb{\overline \omega_{\la}}

\def\si{\sigma}
\def\na(#1){\nabla_{#1}}

\def\evl{\mbox{\em{\rm ev}}}

\def\ev{\mbox{\rm ev}}

\def\im{\mbox{\rm Im}(\cm)}
\def\wt#1{\widetilde{#1}}

\def\ov#1{\overline{#1}}

\def\ga{\gamma_{\varepsilon}}

\def\pls{{\Bbbd PSL}(2,\cx)}
\def\ker{{\rm Ker D}}
\def\cok{{\rm Coker D}}

\begin{document}   
\begin{center}{\bf GENUS ONE ENUMERATIVE INVARIANTS IN $\P^n$\\ 
WITH FIXED $j$ INVARIANT}
\end{center}
\bigskip

\centerline{Eleny  Ionel}
\bigskip

\begin{abstract}\textwidth.4in
We prove recursive formulas for $\tau_d$, the number of 
degree $d$ elliptic curves with fixed $j$ invariant in $\P^n$.  We use
analysis to relate the classical invariant $\tau_d$ to the genus one 
perturbed invariant $RT_{1,d}$ defined recently by Ruan and Tian (the
later invariant can be computed inductively). By considering a sequence
of perturbations converging to zero, we then apply  Taubes' Obstruction Bundle
method to compute the difference between the two invariants. 
\end{abstract}
\bigskip

\setcounter{section}{-1}
\section{Introduction.} 
\setcounter{equation}{0}

A classical problem in enumerative algebraic geometry is to compute the
number of degree $d$, genus $g$ holomorphic curves in $\P^n$ that pass 
through a certain number of constraints (points, lines, etc). 

Let  $\si_d$ denote the number of degree $d$ rational
curves ($g=0$) through appropriate constraints.  For example 
$\si_1(pt,pt)=1$ (since 2 points determine a line). The first nontrivial
cases were computed around  1875 when Schubert,
Halphen, Chasles et al. found  $\si_2$ for $\P^2$ and  $\P^3$. 
Later, more low degree examples were computed in $\P^2$ and $\P^3$,  but the
progress was slow.
  Then in 1993 Kontsevich \cite{K} predicted,
based on ideas of Witten,  that the number $\si_d$ of degree $d$ rational
curves  in $\P^2$ through $3d-1$ points satisfies the
following recursive relation:
\[ \si_d=\ma(\sum)_{d_1+d_2=d} \left[ {3d-1\choose 3d_1-1}d_1^2d_2^2-
{3d-1\choose 3d_1-2} d_1^3d_2\right] \si_{d_1} \si_{d_2} \]
where $d_i\ne 0$, and $\si_1=1$. Ruan-Tian (\cite{rt}, 1994) extended these
formulas  
 for $\si_d$ in any $\P^n$.  
\medskip

When genus $g=1$, the classical problem splits into two totally different
problems: one can count (i) elliptic curves with a fixed complex structure, 
or (ii) elliptic curves with unspecified complex structure (each satisfying 
the appropriate number of constraints). This paper gives recursive 
formulas which completely solve the first of these. 
\medskip

Thus our goal is to compute the number $\tau_d$  of degree $d$ elliptic
curves  in $\P^n$ with fixed $j$ invariant. Classically, the progress on
this problem has been even slower than on the genus one case. Recently,  
Pandharipande \cite{rp} found recursive formulas for $\tau_d$ for the 2
dimensional projective space $\P^2$ using the Kontsevich moduli space of
stable curves. 
 \bigskip

 We will approach the problem from a different direction, using analysis.
Our approach is based on the ideas introduced by Gromov to study symplectic
topology. If  $(\Si,j)$ is 
a fixed Riemann surface, let  
\[\{\mspa \; f:\Si\ra \P^n\; |\; \del f=0, \; [f]=d\cdot l\in
H_2(\P^n,\Z) \;\}/Aut(\Si,j) \]
be the moduli space of  degree $d$ holomorphic maps $f:\Si\ra \P^n$, 
modulo the automorphisms of $(\Si,j)$. Each constraint, such as the requirement
that the image of $f$ passes through a specified point, defines a subset of
this moduli space.
\smallskip

Imposing enough constraints gives a 0-dimensional ``cutdown moduli" space
$\cm_d$. To see whether or not it consists of {\em finitely
many} points, one looks at its  bubble tree compactification $\ov
\cm_d$ \cite{pw}. If the constraints  cut transversely, then all the
boundary strata of $\ov\cm_d$ are at least codimension 1, and thus
empty. Unfortunately, transversality fails at multiply-covered maps or
at constant maps (called {\em ghosts}), so $\ov\cm_d$ is not a manifold. 

 This was a real problem until 1994, when Ruan and Tian considered
the moduli space $\cm_\n$ of solutions of the  perturbed equation:
\[ \del f=\n(x,f(x))\]
and used marked points instead of moding out by $Aut(\Si,j)$. For a generic 
perturbation $\nu$ the moduli space $\cm_\nu$ is 
smooth and compact, so it consists of finitely many points that, counted with
sign, give an invariant $RT_{d,g}$ (independent of $\nu$). 

In $\P^n$,  the genus 0 perturbed invariant
$RT_{d,0}$ is equal to the enumerative invariant $\si_d$.   
The perturbed invariants satisfy a degeneration formula that gives not 
only recursive formulas for the
enumerative invariant $\si_d$ in $\P^n$, but also expresses the higher genus
perturbed invariants in terms of the genus zero invariants \cite{rt}. For
convenience, these formulas are included  in the Appendix. 
\medskip

Unfortunately, when $g=1$, the perturbed invariant $RT_{d,1}$ does not
equal the enumerative invariant $\tau_d$. For example, for $d=2$ curves
in  $\P^2$ the Ruan-Tian invariant is $RT_{2,1}=2$ (cf. (\ref{gendeg})),
while $\tau_2=0$ (there are no degree 2 elliptic curves in $\P^2$). Thus
while the Ruan-Tian invariants are readily computable, they differ from
the enumerative invariants $\tau_d$. One should seek a formula for the
difference between the two invariants.  For that, we take the obvious
approach:  \medskip

Start  with the genus 1 perturbed invariant
$RT_{d,g}$ and consider a sequence of generic perturbations $\n\ra 0$. A 
sequence of $(J,\n)$-holomorphic maps converges 
either to a holomorphic torus or to a bubble tree whose base is a 
constant map (ghost base). 
 Proposition \ref{int} shows that the contribution of the
$(J,0)$-holomorphic  tori is a multiple of $\tau_d$. 
\medskip

We show that the only other contribution comes from bubble trees with
ghost base such that the bubble point is equal to the marked point
$x_1\in T^2$. To compute this contribution, we use the Taubes'   
``Obstruction Bundle" method. Proposition \ref{completion} identifies
the moduli space of $(J,\nu)$-holomorphic maps close to a
bubble tree with the zero set of a specific  section of the obstruction
bundle. Studying the leading order term of this section, we are able to
compute the corresponding contribution (Proposition \ref{zeros}). 
Adding both contributions, yields our main analytic result:
\newcounter{genthe}
\newcounter{gensec}
\setcounter{genthe}{\arabic{equation}}
\setcounter{gensec}{\arabic{section}}
\begin{theorem}\label{gen} 
Consider the genus 1 enumerative invariant  $\tau_d(\b_1,\dots,\b_k)$ in
$\P^n$.
Let $\cu_d$ be the  $n-1$ dimensional moduli space of 1-marked rational 
curves of degree $d$ in $\P^n$ passing through  $\b_1,\dots,\b_k$. 
Let $L\ra \cu_d$ be the relative tangent sheaf, and denote by  
$\wt L\ra \cut_d$ its blow up as in Definition \ref{ltilde}. Then:
\[n_j\tau_d(\b_1,\dots,\b_k)=RT_{d,1}(\b_1\;|\;\b_2,\dots,\b_k)-
\ma(\sum)_{i=0}^{n-1}{n+1\choose i+2} \ev^*(H^{n-i-1})c_1^i(\wt L^*)\] 
where  $H^i$ is a codimension $i$  hyperplane in $\P^n$, $\ev:\cu_d\ra
\P^n$ is the evaluation map corresponding to the special marked point 
and $n_j=Aut_{x}(j)$ is the
order of the group of automorphisms of the complex structure $j$ that fix a
point.   
\end{theorem} 
Theorem \ref{gen} becomes completely explicit provided we can compute 
the top power intersections $\ev^*(H^{n-i-1})c_1^i(\wt L^*)$. We do this
in the second part of the paper, in several steps. For simplicity of 
notation, let 
 \bear
&& x=c_1(L^*) \in H^2(\cu_d,\Z),\quad \wt x=c_1(\wt L^*) \in
H^2(\wt\cu_d,\Z)\\
&&y=\ev^*(H),\quad y\in  H^2(\cu_d,\Z)\quad\mbox{or}\quad 
y\in  H^2(\wt\cu_d,\Z)
\eear
depending on the context. In this notation, Theorem \ref{gen} combined with
(\ref{gendeg})
becomes:
\bear\label{conden}
n_j\tau_d(\; \cdot \;)=
\ma(\sum)_{i_1+i_2=n} \si_d(H^{i_1},H^{i_2},\; \cdot \;)+
\ma(\sum)_{i=0}^{n-1}{n+1\choose i+2}\;\wt x^i y^{n-1-i}\cdot[\wt\cu_d]
\eear
Proposition \ref{pwtc1}  explains how to get 
  recursive formulas relating $\wt x^i y^j$ to  $x^k y^l$ and 
Proposition \ref{pc1} gives recursive formulas for $x^i y^j$ in terms
 of the enumerative invariant $\si_d$. Finally,
the  recursive formulas for $\si_d$ are known (see \cite{rt}, \cite{K}),
so the right hand side of (\ref{conden}) can be recursively computed. 
\bigskip

In the end, we give  applications of these formulas. We
explicitly work out the formulas expressing the number of degree $d$
elliptic curves passing through generic constraints in $\P^2$ and
$\P^3$ in terms of the rational enumerative invariant $\si_d$. For example:

\begin{prop}\label{p^3} For $j\ne 0,1728$, the number $\tau_d=\tau_d(p^a,l^b)$
of elliptic curves in
$\P^3$ with fixed $j$ invariant and passing through $a$ points and $b$
lines (such that $2a+b=4d-1$) is given by:
\bear
\tau_d( \cdot )={(d-1)(d-2)\over d}\si_d(l,\; \cdot )-{1\over
d}\ma(\sum)_{d_1+d_2=d}
d_2(2d_1d_2-d)\si_{d_1}(l,\; \cdot )\si_{d_2}( \cdot )
\eear
where $\si_d(l,\; \cdot )=\si_d(l,p^a,l^b)$ is the number of degree $d$
rational 
curves in $\P^3$
passing through same conditions as $\tau_d$ plus  one more line. The sum 
above is over all decompositions into a degree $d_1$ and a degree 
$d_2$ component, $d_i\ne 0$, and all possible ways of distributing the
constraints 
$p^a,\;l^b$ on the two components.   
\end{prop}
Using a computer program, one then computes specific invariants: for
example, the number of degree 10 tori in $\P^3$ with fixed $j$
invariant and passing through 39 lines is:
\[ 6\cdot 386805671822029784844530703900638969856 \]
when $j\ne 0,1728$. To get $\tau_d$ for $j=0$ or $j=1728$ one simply 
divides the $\tau_d$ computed for a generic $j$ by 3 or 2 respectively. 
\medskip

\non{\bf Acknowledgements.} 
I would like to thank my advisor Prof. Thomas Parker
 for introducing me to the subject and for the countless hours of 
discussions.

\section{Analysis}
\setcounter{subsection}{0} 
\setcounter{equation}{0}
\setcounter{theorem}{0} 
\subsection{Setup}\label{Setup}

Let  $\tau_d$ be the 
 genus one degree $d$ {\em enumerative invariant} (with fixed $j$
invariant)  and $\si_d$ be the 
genus zero degree $d$ enumerative invariant in $\P^n$. Using analytic 
methods, we will compute $\tau_d$ by relating it to the perturbed invariant
$RT_{d,g}$ introduced  by Ruan and Tian \cite{rt}. The later is  defined
as  follows. 
\smallskip

Let $(\Si,j)$ be a genus $g$ Riemann
surface with a fixed complex structure and  $\nu$ an inhomogenous term.  
A $(J,\n)$-{\em holomorphic} map is a solution $f:\Si\ra\P^n$  of   
the equation
\bear\label{jnu}
\del f(x)=\n(x,f(x)).
\eear 
For $2g+l\ge 3$, let $x_1,\dots,x_l$ be fixed marked points on $\Si$, and  
$\al_1,\dots,\al_l $, $\b_1,\dots,\b_k$ be various codimension submanifolds
in $\P^n$, such that 
 \[{\rm  index} \; \del= 
 (n+1)d-n(g-1)=\ma(\sum)_{i=1}^l (n-|\al_i|)+\ma(\sum)_{i=1}^k (n-1-|\b_i|)\]
 For a generic $\n$, the invariant
 \[RT_{d,g}(\al_1,\dots,\al_l\;|\;\b_1,\dots,\b_k)\]
counts the number of $(J,\n)$-holomorphic  degree $d$ maps $f:\Si\ra
\P^n$  that pass through
$\b_1,\dots,\b_k$ with $f(x_i)\in \al_i$ for $i=1,\dots,l$ (for more
details see \cite{rt}).
\medskip

The first part of this paper is devoted to the proof of
Theorem \ref{gen}.
\medskip

\non{\large \bf Outline of the Proof of Theorem \ref{gen}.} 
The proof  is done in several steps. The basic  idea is to start  with the
genus 1 perturbed invariant 
\bear\label{cond}
RT_{d,1}(\b_1\;|\;\b_2,\dots,\b_l)
\eear
and take a sequence of generic perturbations $\n\ra 0$. Denote by
$\cm_{d,1,t\n}$ the moduli space of $(J,t\n)$-holomorphic maps satisfying 
the constraints in (\ref{cond}), and let 
\bear
\cm^\n=\ma(\bigcup)_{t\ge 0}\cm_{d,1,t\nu}.
\eear
As $t\ra 0$, a  sequence of $(J,t\n)$-holomorphic maps converges to a  
$(J,0)$-holomorphic torus or to a bubble tree (\cite{pw}). Let 
$\ov \cm^\n$ denote the bubble tree compactification of $\cm^\n$ (for details
on bubble tree compactifications, see \cite{par}). 
\medskip

 Proposition \ref{int} shows that the number of $(J,t\n)$-holomorphic maps
 converging to a $J$-holomorphic torus is equal to
 \[n_j  \tau_d(\b_1,\dots,\b_k)\] 
where $n_j=|Aut_{x_1}(j)|$ is the order of the group of automorphisms of
the complex structure $j$ that fix the point $x_1$. Namely,
\bear
n_j=\left\{\begin{array}{ll}
 2 &\mbox{ if } j\ne 0,1728\\
6&\mbox{ if } j= 0\\
4&\mbox{ if } j=1728
\end{array}\right.\hskip1in
\eear
These multiplicities occur because  if $f$ is a $J$-holomorphic map, then so 
 is $f\circ \phi$ for any $\phi\in Aut_{x_1}(j)$, but they get perturbed to
different $(J,t\n)$-holomorphic maps.  
\medskip

As $t\ra 0$, there are also  a certain number of solutions converging to 
bubble trees. Because the moduli space of $(J,0)$-holomorphic tori
passing through $\b_1,\dots,\b_k$ is
0 dimensional, the only bubble trees which occur have a multiply-covered
or a ghost base  (for these transversality fails, so dimensions jump up). 
\smallskip

A careful 
dimension count  shows that the  multiply-covered base strata 
are still codimension at least one for  genus $g=1$ maps  in $\P^n$.
(This is {\em not} true for $g\ge 2$.) 
But at a ghost base bubble tree the dimension jumps up by $n$ so 
these strata are $n-1$ dimensional. There are actually  2 such pieces,
corresponding to bubble tree where
 (i) the bubble point is at the marked point $x_1$ and (ii) the bubble point is
somewhere else. To make this precise, a digression is necessary to set up some
notation.  

Let 
\bear\label{cmo}
\cmo_{d}=\{(f,y_1,\dots,y_k)\;|\;  f:S^2\ra \P^n
\mbox{ degree d holomorphic, } f(y_j)\in \b_j \}  
\eear
 be the moduli space of bubble maps, and $\cm_d=\cmo_{d}/G$ 
be the corresponding moduli space of curves, where  $G=\pls$. 
Introduce one special marked point $y\in S^2$  and let   
\bear\label{cu}
\cu_{d}=\{\; [f,y,y_1,\dots,y_k]\;\;|\;\; [f,y_1,\dots,y_k]\in\cm_d \}  
\eear 
 be the moduli space of {\em 1-marked curves} and 
\bear\label{ev}
 \ev:\cu_d\ra \P^n,\;\;\;\ev([f,y,y_1,\dots,y_k]) =f(y).
\eear
be the corresponding  evaluation map. We will use  $f(y)$  to record the image
of the ghost base
\smallskip

 For  generic constraints $\b_1,\dots,\b_k$
 the bubble tree compactification of $\cu_d$ is a {\em smooth }
ma\-ni\-fold that
comes with a  {\em natural stratification}, depending on
the possible splittings into bubble trees and how the degree $d$ and the
constraints $\b_1,\dots,\b_k$  distribute on each bubble. 
\medskip

With this, the  two ``pieces" of the  boundary of $\ov\cm^\n$ 
 are: 
\bear\label{2strata}
  \{x_1\}\ti\ov\cu_d\hskip.2in  \mbox{ and}\hskip.2in T^2\ti \ev^*(\b_1)
\eear

The first factor records the
bubble point, while the image of the ghost base is encoded in the second
factor. For generic constraints, each piece, as well
as their intersection, is a smooth manifold, again stratified. 
\medskip

To see which bubble trees with ghost base appear as a limit of
perturbed tori, we use the Taubes' Obstruction Bundle. This construction must
be performed on the link of each strata. We do this first on the top statum 
of  $\{x_1\}\ti\ov\cu_d$, which consists of bubble trees with ghost base and a
single  bubble. 

In  Section \ref{Approx} we construct  a set of approximate maps 
by  gluing in the bubble. The ``gluing data" $[f,y,v]$ at a 1-marked
curve $[f,y]$ consists of a nonvanishing  vector $v$ tangent to the bubble at
the bubble
point $y$.  Proposition \ref{lg2} shows that the obstruction bundle is  then
diffeomorphic to
$\ev^*(T\P^n)$.
\medskip 

In Section \ref{Gluing} we try to correct the approximate maps to make
them $(J,t\nu)$-holomorphic by
pushing them  in a direction normal to the kernel of the linearized
equation. Those approximate maps that can be corrected to solutions of
the equation (\ref{jnu}) are then identified  with the zero set of a
section $\psi_t$ of the obstruction bundle.   Proposition
\ref{completion} shows that actually all the solutions  of the equation
(\ref{jnu}) are obtained this way, i.e. the end of the moduli  space of
$(J,t\nu)$-holomorphic maps is diffeomorphic to the zero set of the
section $\psi_t$.
\medskip
 
One might be tempted now to belive that the difference between the two
invariants is simply the euler class of the obstruction bundle. But in
fact,  even in generic conditions, the section $\psi_t$ is {\em not a
generic} section of the obstruction bundle. We will see 
that the obstruction bundle has a nowhere vanishing section, so it has a
trivial euler class, while there are examples in which the difference
term is certainly not zero. 
\medskip

Now, to understand the zero set of $\psi_t$ it is enough to look at the 
leading order term of its expansion as $t\ra 0$. By Proposition \ref{psi}
this has the form  $ df_y(v)+t\bar\n$
where $\bar\n$ is the projection of $\n$ on the obstruction bundle. 
\medskip

The construction described above extends naturally to all the
other boundary strata. Each bubble $[f_i,y_i]$ comes with ``gluing data "
$[f_i,y_i,v_i]$, consisting of a vector $v_i$ tangent to the bubble at the
bubble point $y_i$.
But the leading order term of the
section $\psi_t$ depends only on the vectors tangent to the {\em first level
of nontrivial bubbles}. 

More precisely, let  $\cz_h\subset \ov\cu_d$ denote the 
collection of bubble trees for which  the image $u=f(y)$ of the ghost base
 lies on $h$ nontrivial bubbles. Geometrically, the image of a bubble tree
in $\cz_h$ has $h$ components $C_1,\dots,C_h$ that meet at $u$. 
Let  $\cw|_{\cz_h}\ra \cz_h$ be the bundle whose fiber is 
$T_uC_1\oplus\cdots \oplus T_uC_h$.
The  leading order term of $\psi_t$ on $\cz_h$ is a section of  $\cw$, equal to
  
\[a(f,y,v])+t\bar\n\ma(=)^{def}
df_1(y_1)(v_1)+\dots+df_h(y_h)(v_h)+t\bar\n\]
where $([f_i,y_i,v_i])_{i=1}^h$ is the gluing data corresponding to the bubbles
$C_i$, $i=1,\dots,h$.

Unfortunately $\cw\ra\ov\cu_d$ is not a vector bundle (its rank is not
constant).  But if we
blow up each strata $\cz_h$ starting with the bottom one, then the total space
of $\cw$ is the same as the total space of $\wt L$, the blow-up of the
relative tangent sheaf $L\ra \cu_d$. The leading order term of $\psi_t$
descends as a
map $a+t\bar\nu:\wt L\ra \ev^*(T\P^n)$. Moreover, 
$\bar \nu$ doesn't vanish on $\im=ev_*(\ov\cu_d)$ so it induces
a splitting on the restriction 
\[T\P^n/\im=\cx \lg \bar\n\rg\oplus E.\] 
 Finally, we put  all these pieces together in Proposition \ref{zeros} to
prove that the number of $(J,\n)$-holomorphic maps converging as $\nu\ra
0$ to the boundary strata $\{x_1\}\ti\ov\cu_d$ is given by the Euler class 
$c_{n-1}(\ev^*(E)\otimes\wt L^*)$.
\medskip

In Section \ref{Other} we show that the other boundary strata $T^2\ti
\ev^*(\b_1)$ 
gives trivial contribution, concluding the proof of the Theorem \ref{gen}.

\subsection{The Approximate gluing map}\label{Approx}

Let $\cu_d$ be the moduli space  of 1-marked  rational
curves of  degree $d$ passing through the conditions
$\b_1,\dots,\b_k$. 
 In this section we  construct a set of approximate maps starting
from $ \{x_1\}\ti\ov\cu_d $, the first boundary strata in (\ref{2strata}).
We will use a:
\step{Cutoff function} Fix a smooth cutoff function
$\b$ such that $\b(r) =0$ for
$r\le 1$ and  $\b(r) =1$ for $r\ge 2$. Let  $\b_\la (r) =
\b(r/\sla )$. Then $\b_\la$ has the following properties:
\best  
| \b_\la | \le 1\; , \;\;| d \b_\la | \le 2/\sla\;\;  \mbox{ and }\;\;
d\b_\la \;\;\mbox{ is supported in } \sla \le r \le 2 \sla  
\eest 
\step{The definition of the approximate gluing map on the top stratum} Let
$\cn$ denote the top stratum of $\{x_1\}\ti\ov\cu_d$.  
 First we need to  choose a canonical representative of each bubble
curve $[f,y]\in\cn$ (recall that $f(y)$ is the image of the ghost base).  
Using the $G=\pls$ action, we can assume that $y$ is the North pole and 
 $f$ is centered on the vertical axis, which  leaves a $\cx^*\cong
S^1\ti\R_+$ indeterminancy. To break it off, include as  gluing data a 
unit vector tangent
to the domain $S^2$ of the bubble at the bubble point $y$. The frame
bundle 
\bear\label{fr}
 Fr=\{\; [f,y,u] \;|\; [f,y]\in\cu_d,\; u\in T_yS^2,\;|u|=1\}
\eear 
models the link of $\cn$. The notation $[f,y,u]$ means the
equivalence class under the action of $G$ given by:
\best
g\cdot(f,y,u)=(f\circ g^{-1}, \;g(y),\; g(u) )
\eest
where the compact piece $SO(3)\subset G$ acts on the unit frame $u$ 
by rotations and the noncompact part acts trivially. 
\smallskip

Fix a nonzero vector $u_1$ tangent to the torus at $x_1$. 
This determines 
an identification $T_{x_1}(T^2) \cong \cx $ such that $u_1=1$, giving
local coordinates on the torus at $x_1=0$. Similarly, let $u_0$ be a
unit vector tangent to the sphere $S^2$ at the north pole and consider
the identification 
\bear\label{tt2} (T_{x_1}T^2, u_1)\cong (T_NS^2,u_0)\eear
that induces natural coordinates on the sphere via the stereographical
projection (such that $N=0,\;u_0=1$). These choices of local coordinates on the
domain of the bubble tree will be used for the rest of the paper.
Fix also a metric on $\P^n$ such that we can use normal coordinates up to 
radius 1.


To glue, one needs to make sure that only a small part of the
energy of $f$ is concentrated in a neighbourhood of $y$. 
The convention in \cite{pw} is to  rescale $f$ until $\ep_0$
of its energy is distributed in $H_y$, the hemisphere centered at $y$. 

But since the constructions in the next couple of sections involve quite a few
estimates, we prefer to do a different rescaling, that will simplify the
analysis.  Choose a representative of $[f,y,v]$ such that 
\bear\label{uniquerep1}
 y=0, \;\; u=1,\;\; f \mbox{ centered on the vertical axis }
\eear
Since on the top strata $[f,y]$ cannot be a ghost, such
representative  is uniquely determined up to a rescaling factor
$r\in\R_+$. We will choose this rescaling factor such that moreover
\bear\label{d2f}
{\rm max}\{\;|\nabla^2f(z)|,\; |z|\le 1\}\le 2
\eear 
Note that if the degree of $f$ is not 1, then imposing the extra condition 
\bear\label{uniquerep2}
{\rm max}\{\;|\nabla^2f(z)|,\; |z|\le 1\}=2
\eear 
determines uniquely the representative. To see this, choose some  
representative $\wt f$  as in  (\ref{uniquerep1}) and look for a
map $f(z)=\wt f(rz)$ satisfying also (\ref{uniquerep2}). The uniqueness
comes from the fact that the map  
$s(r)={\rm max}\{\;|\nabla^2\wt f(z)|,\; |z|\le r\}-{2/r^2}$ is decreasing. 
\smallskip

If the degree of $f$ is 1, (i.e. the image curve is a line), then we could
 replace (\ref{uniquerep2}) by say $ |df(0)|=1 $ and still have 
(\ref{d2f}) satisfied. 
\medskip

Finally, the {\em approximate gluing map} 
\bear\nonumber
&&\ga :Fr\ti (0,\ep) \ra \cm aps(T^2,X)  \\
&&\ga(\;[f,y,u],\;\la)=\fl
\eear
 is constructed  as follows: Choose the unique representative of $[f,y,u]$ 
satisfying (\ref{uniquerep1}) and (\ref{uniquerep2}). The  approximate map 
$\fl$ is obtained by gluing
to the constant map   $f(y)$ defined on $T^2$  the bubble map $f$ 
rescaled by a factor of $\la$ inside a disk $D(0,\sla)\subset T^2$,
\[ \fl (z) = \b_ \la(|z|) f  \l( {\la \over z} \r) \]
where the multiplication is done in normal coordinates at $f(0)$.  
 Let $Gl=Fr\ti (0,\ep)$ denote the set of gluing data.



\step{Weighted Norms} On the domain of $f_{\la}$ we will use the
rescaled metric $g_\la=\tel^{-2}dzd\bar z$, where
\[   \tel(z)= (1-\b_\la(z)\;)(\la + \la^{-1} |z|^2)+\b_\la(z)\] 
Define 
\best
  \| \xi \lpi&=&  \left( \int |\xi|^p \tel^{-2}  +  
 | \nabla \xi |^p \tel^{p-2} \right)^{1/p} \mbox{ for } \xi 
\mbox{ vector field along } f_{\la} \mbox{ and}\\
  \| \eta \lp&=&\l( \int   |\eta |^p  \tel ^{p-2} \r)^{1/p}
 \mbox{ for } \eta \mbox{ 1-form along } f_{\la} 
\eest

\non The weighted norm of a vector field or 1-form on $f_{\la}$ equals 
its  usual norm off $B(0,2\sla)$ and on $ B(0,\sla)$ it is equal with the
norm of its  pulled back on $S^2$ via a rescaling of factor $\la$. The
usual Sobolev embeddings hold for this weighted norms with constants
independent of $\la$. 

\begin{lemma}{\label{lg1o}}
There exists  $\ep_{0} > 0$ and  constants $C > 0 $  such that for any 
$p \ge 1$ and  $\la  \; \le \;   \ep_{0}$:  
 \be\label{fla}
\| df_{\la} \lp  \; \le \;\;\;   C \mbox{ and }\;\;\;
\| \del f_{\la}  \lp  \; \le \;   C \la^ {1/p}
\end{equation}
\non Moreover on the annulus A: $\{\sla  \; \le \;   |z|  \; \le \;  2\sla\}$
we have the following  expansion:
\be{\label{pe1}}
\del \fl = {\sla \over |z|}\; d\b \cdot d f (y)(u) + \o (\la) 
\end{equation} 
The estimates are uniform on $Gl\ra \cn$.  
\end{lemma}

\pf Let $B$ be the disk $|z| \le \sla$. Note that $d \fl$ vanishes for $|z| \ge
2\sla$  and by the definition of the weighted norm on $B$,
\[ \|d \fl \|_{p,\la,B} = \|d f \|_{p,B}\]
But (\ref{d2f}) implies that 
\bear\label{df}
{\rm max}\{\;|df(z)|,\; |z|\le 1\}\le 2 
\eear
 In the same time, $ \del \fl=0$ outside A.
Hence we need only to consider what happens in A. But on $A$  
 
\[ |\del \fl |  \; \le \;   C | d \fl | \; \le \; 
 C( |d \bl|\;| f | + |\bl |\;| df |)\;
 {\la\over |z|^2}   \; \le \;    
C{1 \over \sla}  \ma( \sup )_{B}|f | +C  \; \le \;   C   \]

\non since $\ma( \sup )_{B}|f (z)|  \; \le \;   
\sla\; \ma( \sup )_{B}| d f |   \; \le \;  2 \sla\;$ 
in normal coordinates on $\P^n$ at $f(y)$. This concludes the first part of
the proof. For the second part, notice that on $A$
\[ \del \fl = \del \bl \cdot f  + \bl \cdot \del f  =
 {1\over \sla}\; d \b \;{z \over |z|}\cdot f  \l( {\la \over  z}\r) \]
\non since  $f $ is holomorphic. But using (\ref{d2f}) in normal
coordinates on $\P^n$ at $f(y)$ and $y=0$, we get 
$|f(z)-f(0)-df(0)(z)|\le 2|z|^2$  so
\[ f \l({\la\over z}\r)={\la\over z}\cdot df_y(u)+\o(\la)\;\;\mbox{  on }\;A \]

Substituting this in the formula for $\del \fl$ we obtain (\ref{pe1}).\qed
\medskip

\step{Extending the approximate gluing map}
 The approximate gluing map extends
naturally to the bubble tree compactification $\ov\cu_{d}$ of the
moduli space of 1-marked  curves. For simplicity, let $\cn$ denote some
 boundary stratum modeled on a bubble tree $B$ and  corresponding to a 
certain  distribution of the degree $d=d_1+\dots+d_m$ on the bubbles. If 
  $[f_i,y_i],\;i=1,\dots,m$ are the bubble curves corresponding to the bubble
map $f:B\ra\P^n$, then the gluing data $Gl$ 
is a collection of unit vectors tangent to each sphere in the domain at the
 corresponding bubble point together with gluing parameters: 
\bear\label{gl}
Gl=\{\; (\;[f_i,y_i,u_i],\la_i\;)_{i=1}^m\;|\; u_i\in T_{y_i}S^2, 
\;|u_i|\ne 0, \la_i\le \ep\}
\eear

Note that as long as $f_i$ is not a constant map, then we can choose a
unique reresentative of $[f_i,y_i,u_i]$ as in (\ref{uniquerep1}),
(\ref{uniquerep2}). Then Lemma \ref{lg1o} extends naturally to $\cn$ to give

\begin{lemma}\label{lg1} With the notations above, let $\fl$ be an approximate

gluing map, and  $A_1,\dots A_m$ be the corresponding annuli of radii 
$\la_i$ in which the cutoff functions are supported. Then for $\ep$
small enough, there exists a constant $C$ such that:
\[\| df_{\la} \lp  \; \le \;\;\;   C\;,\;\;\;\;
\| \del f_{\la}  \lp  \; \le \;   C \la^ {1/p} \] 
Moreover,   $\del\fl=0$
except on the annuli $A_i$ that correspond to nontrivial bubbles, where
\bear\label{aprglbd}
\del \fl=
-{\sla_i\over |z|}\; d\b \cdot d f_i (y_i)(u_i) + \o (\la_i) 
\eear 
The estimates above are uniform on $Gl\ra \cn$. 
\end{lemma}

We will see later that most of the important information is
encoded in the first level of nontrivial bubbles.

\subsection{The Obstruction Bundle}\label{Obstr}

In order to see which of the approximate maps can be corrected  to
solutions of the equation $\del f=\nu$  we  need first  
to understand the behaviour of the linearization of this equation over
the space of approximate solutions.  
\smallskip

Recall that transversality fails at a bubble tree with ghost base, so
the linearization at such bubble tree is not onto. The cause of that is
the ghost base. Thus we  start by analysing the ghost maps:
 
Consider the moduli space of 
holomorphic maps $f:T^2\ra \P^n$ representing  $0\in H_2(\P^n )$. 
Obviously, the only such maps are the constant ones (ghosts). 
If  $D_u$ is the linearization of the section 
$ \del : \cm aps(T^2,\P^n ) \ra \La $
 at $f:T^2\ra \P^n $, $f(x)=u$ a constant map, then 
\best
 \mbox{ index D}_u&=&\mbox{dim }\ker_u- \mbox{dim }\cok_u =
 c_1(0) + n(1-1)=0
\eest
 and 
\[ \cok_u = \ho(f)\cong T_u\P^n\;\;\;\mbox{ (canonically) }\]
\non since $f^*(T\P^n)$ is a trivial bundle, so the elements 
$\om\in \ho(f)$ are constant on the torus, i.e. have the form  $\om=Xdz$ for 
some $X\in T_u\P^n$. 
\bigskip

Now if $f:B\ra \P^n$ is a bubble tree map whose base is a ghost torus
$u=f(y)\in \P^n$, let  $D_f$ be the linearization at $f$ of the section 
$ \del : \cm aps(B,\P^n ) \ra \La$. Then
\best
 \mbox{ index D}_f&=&\mbox{dim }\ker_f- \mbox{dim }\cok_f =-1
\eest
To describe $\cok_f$ we will use the following:
\begin{defn} \label{obstrbd} If $f:B\ra \P^n$ is as above, let

$B_1\subset B$ consist of the domains of all the ghost bubbles with image
$f(y)$,
 
$B_2=B -B_1$ and 
 
$\wt B\subset B$  denote the first level of bubbles that are not in $B_1$. 
\end{defn} 
Then  $\cok_f $ is $n$ dimensional, consisting of 1-forms $\om$ such that 
\[ \om=\left\{ \begin{array}{ll}
X dz &\mbox{ on } B_1\\
0&\mbox{ on } B_2
\end{array}\right.\]
for some $X\in T_u\P^n$. In particular, there is a natural isomorphism 
\bear\label{obB}
\begin{array}{c}
\cok \;\;\cong \;\;\ev^*(TP^n)\\
\searrow\;\;\;\;\;\;\;\swarrow\;\;\\
\ov\cu_d
\end{array}
\eear
where $\ev:\ov\cu_d \ra \P^n$ is the evaluation map. 
Since the moduli space of bubble trees $\ov{\cu_d}$ is compact, there
exists a constant $E>0$ such that $D_f D^*_f$ has  a  zero eigenvalue with
multiplicity $n$, and all the other eigenvalues are greater than $2E$
for all $f\in \ov \cu_d$.
\bigskip

When $\fl$ is an approximate map, let $\Dl$ be the 
linearization of \[\del : \cm aps(T^2,\P^n ) \ra \La \] at $\fl$  
and  $ \Dls$ its $L^2$-adjoint with respect to the metric 
$g_\la$ on $T^2$. Then $\Dl$ is not uniformly invertible. More
precisely, 
\begin{lemma}\label{lg2}
For $\la > 0$ small, the operator 
$\Delta_{\la}= \Dl  \Dls $ has exactly $n$ eigenvalues
of  order $\sla$ and all the others are greater than E. Moreover, over
the  set of gluing data $Gl$, the span of low eigenvalues
 \best
\begin{array}{ccc}
 \lal(f_\la) &\hookrightarrow& \Lal\\
&&\downarrow   \\
&&Gl\end{array}
\eest  
is a $n$-dimensional vector bundle 
(called the {\rm Taubes obstruction bundle}),  
naturally isomorphic to the bundle
\[ \evl^* (T\P^n )\ra Gl\] 
where $\evl:  Gl \ra\P^n $ is the evaluation map. 
\end{lemma}
\pf The proof is more or less the same as the one Taubes used for  the 
similar result in the context of Donaldson theory, \cite{t2}. 
 For each gluing data in $Gl$,  by cutting and pasting eigenvectors  
we show that the eigenvalues of $\Delta_\la=\Dl \Dls$ are 
$\o (\sla)$ close to those of $\Delta_u=D_u D^*_u$, 
where $u$ is the point map in the base of the bubble tree.

Take for example the top stratum of $\ov\cu_d$. 
Choose $\{\om_i,\; i=1,n\}$ a local orthonormal 
base of $\cok\cong\ev^*(T\P^n)$  and define
\bear\label{taubes}
 \oml^i (z) = \b\l({z\over 2\sla}\r) \om^i(z) 
\eear
 A straightforward computation shows that:
\bear
\| \Dls  \om_\la \|_{2,\la}   &\le &   \la^{1/4} \|\om_\la
\|_{2,\la}\\ 
\lg \omb^i,\omb^j \rg _{2,\la}&=&\de_{ij}+{\o}(\la)
\eear  
The Gramm-Schmidt orthonormalization procedure then provides $n$ 
eigenvectors $\omb^i$ for $\Delta_{\la}$ with eigenvalues $\o (\sla)$
such that  
\[   \omb^i =\oml^i + \o(\la) \]
The construction above  extends
naturally to the other substrata of $\ov\cu_d$. Note that for example 
 when $B_1$ has other  components besides $T^2$ then $\omb$ is equal to 
$\om$ not only on the ghost 
base, but on all $B_1$ and is extended with 0 starting from the first level
of nontrivial bubbles.
\medskip

An adaptation of Taubes argument from \cite{t1} 
shows that  there are at most $n$ low  eigenvalues of $\Delta_\la$.  
Therefore there is a well defined splitting 
\[\lla(f_\la) = \lal(f_\la) \oplus \lae(f_\la)\]
 The definition (\ref{taubes}) combined with (\ref{obB}) 
 provides the isomorphism $\Lal \cong \ev^*(T\P^n )$, concluding the proof.
\qed

\step{The partial right inverse of $D_\la$}  The  restriction of $D_{\la}
D^{*}_{\la} $ to $\Lae$ is invertible (since all its eigenvalues are at least
E). Define  $  P_\la$ to be the composition of the $L^2$-othogonal projection
$\La \ra \Lae$ with the operator
 $ D^*_\la ( {D_{\la}  D^{*}_{\la}})^{-1} $ on ${ \Lae}$.
 Then  \bear P_\la : \lla(f_\la)    \ra  \lo(f_\la) \eear 
is  the {\em partial
right inverse of}  $D_\la$ and satisfies the uniform estimate:
 \be   \|P_\la \eta\lpi     \; \le \;   E^{-1} \| \eta \lp 
\end{equation}

\non We will denote  by $\pr^{\fl} : \lla(f_\la) \ra \lal(f_\la)$ the 
projection onto the fiber of the obstruction bundle.

\subsection{The Gluing map}\label{Gluing}

 The next step is to correct the approximate gluing map to take values in 
the moduli space $\cm_{t\n}$ of solutions to the equation
\be\label{eg0}
 \del f(x)=t\cdot \n(x,f(x)) 
\end{equation}
\non  where $\n$ is generic and fixed and  $t$  is a small parameter.

If $f_\la$ is an approximate map, use the exponential map to write any 
nearby map in the form $f=\exp_{\fl} (\xi) $, for some  correction 
$\xi \in \lo(\fl)$. 
Let $\Dl$ be the linearization of the $\del$-section at $\fl $ so
\be\label{del0}
\del f = \del \fl + \Dl(\xi) + Q_\la(\xi) 
\end{equation}
\non where $ Q_\la$ is quadratic in $\xi$. Similarly, 
\[\n(x,f(x))= \n(x,\fl(x)) + d\n(\xi)+\widetilde Q_\la(\xi)\] 
so equation  (\ref{eg0}) can be rewritten as:
\bear\label{eg00}
{ \Dl(\xi)+N_\la(\xi,t)  =  t\n(x,\fl(x)) -\del \fl  }
\eear
where $N_\la(\xi,t)=Q_\la(\xi)-td\n (\xi)-t\widetilde 
 Q_\la(\xi)$ is  quadratic in $(\xi, t)$.
\medskip

\non The kernel of $\Dl$ models the tangent directions to the space of   
approximate maps, so it is natural to look for a correction in the normal 
direction. More precisely, we will consider the  solutions of  
(\ref{eg00}) of the form 
\bear
f=\exp_{f_\la} ( P_\la \eta)\;\;\; \mbox{ where  }\;\;\;\pr(\eta)=0
\eear
Since $ \Dl ( P_\la(\eta))= \eta$ for such $\eta$, then equation 
(\ref{eg00}) becomes
\be\label{eg1}
\eta + N_{\la,t} ( P_\la \eta) = t \n-\del f_\la 
\end{equation}  
The existence of a solution of (\ref{eg1}) is a standard aplication of the 
Banach fixed point theorem combined with the estimates in the previous 
sections.  
\begin{lemma}\label{l3}
There exists a constant $\de >0$ (independent of $\la,t$) such that for 
$t$ small enough and for any $\alpha \in \lla(f_\la) $ so that 
$\| \alpha \lp < \de/2$ the equation:
\[ \eta + N_{\la,t} ( P_\la \eta) = \alpha \]
has a unique small solution $\eta \in \lla(\fl)$ with $\| \eta \lp < \de $.
 Moreover, 
\[\| \eta \lp < 2 \| \alpha \lp \]
 and if $\alpha$ is $C^\infty$, so in $\eta$.
\end{lemma}
\pf  Apply the contraction principle   to the operator 
\[ T_\la : \lla(f_\la) \ra \lla(f_\la)  \]
 \[ T_\la \eta = \alpha -   N_{\la,t} ( P_\la \eta) \]
\non defined on a small ball centered at 0 in the Banach space 
$\lla(\fl)$ with the weighted Sobolev norm $L^p_\la$. 
To prove that T is a contraction we note that: 
\[  \| T_\la \eta_1 - T_\la \eta_2 \lp = 
\| N_{\la,t} ( P_\la \eta_1)- N_{\la,t} ( P_\la \eta _2) \lp \] 
\non and use some estimates of Floer. He proved in \cite{F} that  for 
the quadratic part $Q$ of (\ref{del0}), 
there exists a constant $C$ depending only on 
$\| df \lp$ such that:
\bear\label{floer1}
\| Q_f(\xi_1)- Q_f(\xi_2) \lp   \; \le \;   C \; (\; \|\xi_1 \lpi + \|\xi_2
\lpi ) \| \xi_1-\xi_2 \lpi 
\eear
\bear \label{floer2} 
\| Q_f(\xi)  \lp   \; \le \;   C \; \|\xi \|_{\infty,\la} \cdot \|\xi \lpi
\eear
(Floer's estimates are for the usual Sobolev norm, but the same proof
goes through for the weighted norms.) 
 Since   $\| d\fl \lp$  is uniformly bounded by Lemma $\ref{lg1}$, the same
constant $C$ works for all $ f_\la \in Im(\ga)$. Moreover, for  $t$ very small
the same estimates hold for the nonlinear part $N_{\la,t}$.

Hence by (\ref{floer1}):
\best
 \| T_\la \eta_1 - T_\la \eta_2 \lp & \le& 
C \;(\; \| P_\la \eta_1 \lpi + \| P_\la \eta_2 \lpi) \|  P_\la (\eta_1-\eta_2)
\lpi           \\   
    & \le &  C/E^2 \;( \;\| \eta_1 \lp + \| \eta_2 \lp ) \cdot\| \eta_1-\eta_2
\lp.
\eest
Choosing $\de < E^2/(4C)$ this implies  
\[ \| T_\la \eta_1 - T_\la \eta_2 \lp   \; \le \;   1/2 \;\| \eta_1-\eta_2 \lp
\]
 for any $\eta_1,\eta_2 \in B(0,\de)$. Moreover, since 
 $ \| T_\la (0) \lp \le \de/2 $ then $  T_\la :B(0,\de)   \ra   B(0,\de)$ 
is a contraction.  
Therefore $T_\la$ has a unique fixed point $\eta$ in the ball such that
moreover
\[ \|\eta \lp \;\le \;\|T_\la\eta -T_\la(0)\lp + \|T_\la(0)\lp\;\le 1/2
\;\|\eta \lp+ \| T_\la(0)\lp\] 
\non so  $\|\eta \lp \le 2\;\|T_\la(0)\lp=2\|\alpha\lp$. Elliptic 
regularity  implies that $\eta$ is smooth when $\al$ is.\qed
\begin{cor}\label{defeta}
For $t,\la$ small enough,  equation (\ref{eg1}) has a unique small solution 
$\|\eta\lp   \; \le \;   \de$. Moreover,  
 \[ \|\eta\lp   \; \le \;   C( t|\nu | + \la^{1\over p}). \]
\end{cor}
\pf Follows immediately from Lemmas \ref{lg1} and \ref{lg2}  and the 
estimate 
\[ \| \alpha \lp = \| \;t\nu -\del \fl \lp   \; \le \;   
t|\nu| + C \la^{1\over p}. \hskip.5in  \Box\] 

\step{The gluing map}  Let $ Gl$ be the set of gluing data. The {\em
gluing  map} is defined by
\[ \bar \ga :Gl   \ra   \cm aps(T^2, X) \]
\[  \bar \ga ([f,y,u],\la) = \bar f_\la =\exp_{f_\la}(P_\la \eta ) \]
\non where $\eta=\eta(f,y,u,\la)$ is the unique solution to the equation
(\ref{eg1}) given by Corrolary \ref{defeta}.

 By construction, $  \bar \ga $ is a local diffeomorphism onto its image.
Moreover, if  $\pr^{f_\la}(\eta)=0$ then $\bar \fl$ is actually a solution of 
(\ref{eg0}).

\step{The obstruction to gluing}  The section
\[ \psi_t :Gl \ra  \lal(f_\la) \;\;\;\mbox { given by } \]
\[ \psi_t (f,y,u,\la)  = \pr^{f_\la}(\eta)= \pr ^{f_\la}( t \n - \del f_\la ) -
\pr ^{f_\la}( N_{\la,t} ( P_\la \eta))  \]
will be called   the {\em obstruction to gluing}.
\non Let $Z_t =\psi_t^{-1} (0)$ be the zero set of this section. By applying 
 the gluing construction to bubble trees in $Z_t$ we get a subset
 of the moduli space  $ \cm^{t \n}$.

\subsection{Completion of the construction}\label{Compl}

We have seen in the previous section that applying the gluing construction 
to the bubble trees in the zero set $Z_t$  we will get 
elements of the moduli space $\cm_{d,1,t\n}$. It is not clear yet why all the 
elements of this moduli space close enough to the boundary stratum 
$\cn$ can be  obtained   by  the gluing procedure. The purpose of this section
is to 
clarify this issue. 
\medskip

 Recall the construction of the gluing map: Starting with a bubble tree we 
 glue in the bubble to obtain an approximate map $\fl$. Then we   
correct $\fl$ by pushing it in a direction normal to the kernel of $D_\la$ 
in order  to get an element of the moduli space $\cm^{t\n}$. 
The key fact here is  that the kernel of the linearization  $D_\la$ models the

 tangent space to the approximate maps, and therefore, at least in the 
 linear model, it is enough to look for solutions only in a normal direction.
 For the construction to be complete though, we need to show that the same
 thing is true for the nonlinear problem.
\medskip

 More precisely, we will  show that for $t$ small, all the elements of the
 moduli space $\cm_{d,1,t\n}$ 
 close to the boundary stratum $\cn$ can be reached starting with an 
approximate map and  going out in a normal direction. The proof of the 
following Theorem is an adaptation of the proof for the same kind of 
result in the context of Donaldson theory \cite{dk}. It is pretty 
technical and we include it just for continuity of the presentation. 

\begin{theorem}\label{completion}  The end of the moduli space 
$\cm_{d,1,t\n}$ close to the
boundary strata $\cn$ is diffeomorphic to the zero set of the section
$\psi_t$. More precisely, for  $\de $ and $t $ small enough, there
exists an isomorphism 
\[\cm_{d,1,t\n} \cap \cud  \cong   \psi^{-1}_t(0)\hskip.6in \mbox{where}\]
\[\cud =\{ f:T^2\ra X \mid \exists \fl \;s.t.\; 
f=\exp_\fl(\xi) ,\; \| \xi \lui \le \de,\; 
\| \del f\lu \le  \de^{3/2} \}   \]
and $ f_\la   \in  {\rm Im} \ga$ is some approximate map.
\end{theorem}
\pf The proof consists of 2 steps.  First, Lemma \ref{nbd} shows
that $\cud$ is actually a neighborhood of $\cn$ in the  bubble tree
convergence topology. Second, recall that in  constructing  the section
$\psi_t$ 
we were looking for solutions of the equation  (\ref{eg0}) that have the form 
\be\label{form}
 f=\exp_{\fl}( P_\la \eta )\;\;\mbox{ for some }\; \|\eta \lu  \; \le \; \de 
\end{equation} 
 To prove the Theorem it is enough to show that for $t$ small, all the
solutions 
of the equation (\ref{eg0}) can be written in the form
 (\ref{form}). This is a consequence of Proposition \ref{prop}.
\begin{lemma}\label{nbd}
$\cud\cap \ov{\cm^\n}$ is a neighborhood of $\cn$ in the bubble tree 
convergence topology. More precisely, for any  $(J,t\n)$-holomorphic map f
close to the boundary strata $\cn$ there exists an approximate map $\fl$ such
that f can be written in the form 
\[ f=\exp_\fl( \xi )\;\;\mbox{ for some }\; \|\xi \lui  \; \le \; \de   \] 
\end{lemma}
\pf  By contradiction, assume there exists a sequence  
$f_n$ of $(J,t_n\n)$-holomorphic maps for ${t_n}\ra 0$ such that 
$f_n$ do not have the required property. 
By the bubble tree convergence Theorem (\cite{pw}) there exists a bubble 
tree $f $ such that $f_n\ra f $ uniform on compacts. Moreover, after 
rescaling the functions $f_n$ by some $\la_n$, this becomes a 
$L^{1,2}$-convergence (cf. \cite{pw}). But this is equivalent to saying that
$f_n$ is 
$L^{1,2,\la_n}$ close to $f$. In particular, for $\la$ small enough, 
$f_n$ is  $L^{1,2,\la_n}$ close to $f_{\la_n}$, which gives a
contradiction. \qed

\begin{prop}\label{prop}
For small enough $\de,t$  any map in $\cud$ can be represented in the form 
\[ f=\exp _{f_\la} (P_{\la} \eta)\;\mbox{  for some }\;f_\la \in {\rm Im} \ga,
 \;\| \eta  \lu  <  \de \;\mbox{ and }\; \pr^{\fl} (\eta)=0 \]
\end{prop}
\pf We will use the continuation method. The key fact is that 
 a neighborhood of $f_\la $  in $ {\rm Im} \ga $ is modeled by 
$ \Lol $ and that the image of $P_\la$ spans the normal directions to 
$ {\rm Im} \ga  $. 

Let $f \in \cud$. By definition, there is $f_\la \in {\rm Im} \ga$ such that $
f=\exp_{f_\la }\xi $, where $ \| \xi \lui< \de $. Consider the path 
$f_s=\exp_{f_ \la}( s \xi) $. Let 
 \bear
 S =\{ s \in [0,1] \mid 
\exists f_{\la_s}  \mbox{ and } \|\eta_s \|_{p,\la_s} < \de 
\mbox{ such that }  f_s= \exp _{f_{\la_s}} (P_{\la_s} \eta_s) \}.
\eear

\non  Note that by definition $f =f_\la=\exp_\fl(0)$ so  $0\in S$. 
 We will show that S is both open and closed and since it is nonempty, 
$1 \in S $.
\medskip

\non {\bf S is closed}. The only open condition in the definition of S 
is $\|\eta_s \|_{p,\la_s} < \de $. But since 
\bear\nonumber
 \del f_s &=& \del f_{\la_s} + D_{\la_s}(P_{\la_s}\eta_s ) + 
N_{\la}(P_{\la_s}\eta_s) \hskip.3in \mbox{then}\\ \nonumber
 \eta_s&=& \del f_s - \del f_{\la_s} - N_{\la_s}(P_{\la_s}\eta_s)  
\hskip.5in \mbox{so}\\ \nonumber    
 \| \eta_s  \lu     &\le&       
 \| \del f_s  \lu  + \| \del f_{\la_s}  \lu  +
{C \over E^2} \| \eta_s \lu^2  \\
\label{etas}
   &\le & \| \del f_s  \lu  + C \sla + C \| \eta_s  \lu  ^2 
\eear 
We need to estimate $\|  \del f_s \lu$. But 
 \bear\nonumber 
&&\del f_s = \del f_{\la} + s D_{\la}(\xi ) +N_{\la} (s\xi) 
\hskip.2in\mbox{and}\\ \nonumber
&&\del f_1 = \del f_{\la} + D_{\la}(\xi ) +N_{\la} (\xi) 
\hskip.2in\mbox{so}\\ \label{delfs2}
&&\del f_s = s\del f_{1}+(1-s)\del f_{\la}+N_{\la} (s\xi) - sN_{\la} (\xi)
\eear
\non The relation (\ref{delfs2}) combined with the estimate (\ref{floer2})
gives
 $\;\|N_{\la} (\xi) \lu \le C \;\| \xi \lui^2 $ so
\[ \| \del f_s  \lu  \; \le \;  \| \del f_1  \lu  + \| \del f_\la \lu  +
2\;C \; \|\xi\lui^2  \; \le \; \sla +\de^{3/2} + C\; \de^2\] 

\non Therefore for $\la<<\de$,
\be\label{delfs}
\| \del f_s  \lu  \; \le \; 2\;C\;\de^{3/2}
\end{equation}
 Using (\ref{delfs}) in (\ref{etas}) we get
\[  \| \eta_s  \lu     \; \le \;  2\;C\;\de^{3/2} + 
C \sla + C \| \eta_s  \lu  ^2 \]

\non For small $ \la \; \le \; \de^{ 3}$ , the constraint $ \| \eta_s  \lu   <
\de $ implies
  $ \| \eta_s  \lu   < \de/2  $ so it is a closed condition too. \qed
\medskip

\non {\bf S is open}. Assume that $s_0 \in S$, i.e. there exists an approximate
map  $f_{\la_0}$ such that $f_{s_0}=\exp _{f_{\la_0}}(P_{\la_0}(\eta_0))$. We 
 will show that $s\in S$ for $s$ sufficiently close to $s_0$. For that we need
to find an approximate map $f_{\la_s}$ and an $\eta_s \in \Loe$ such that:
 \be\label{long}
f_s=\exp_{\fl}(s\xi) =\; \exp _{f_{\la_s}} (P_{\la_s} \eta_s)
\end{equation}
It is enough to prove that the linearization of the equation 
(\ref{long})  is onto at $s_0$. First we prove that:

\begin{lemma}\label{neigh}
 A small neighborhood $\cn_\de$ of  $f_ \la $ in ${\rm Im}  \ga$ is modelled by
$  \Lol   $. 
 More precisely, there is a well defined map 
$ g: \Lol  \ra \Lae$ such that any approximate map 
$f  \in {\rm Im}  \ga$  has the form 
$f=\exp _{f_ \la } (\ze+ P_\la g(  \ze))$  for some  $\ze \in \Lol$, $ \| \ze
\lui  \; \le \; \de$.
\end{lemma}
\pf The first statement is an immediate consequence of the way we constructed
the approximate maps.  For the second part, notice that any $f  \in {\rm Im}
\ga $ close to $f_\la$ can be written in the form $f=\exp _{f_ \la } (\chi)$,
with $\chi$ small.  Let $ \chi = \ze + P_ \la  \eta $ be the orthogonal
decomposition of $\chi$ in $\Lol \oplus \Loe $, where 
$\eta \in \Lae$ (recall that $P_\la: \Lae \ra \Loe $ is an isomorphism). Using
the same techniques as in Section \ref{Gluing} we can prove that for any
$\ze\in \Lol$  there exists a unique solution $\eta=g(\ze)$ to the equation
\[  \eta + N_\la(P_\la\eta) = \del f\]
\non which concludes the proof of Lemma. \qed
\bigskip

Since the notations are becoming cumbersome, we will illustrate for simplicity
the case $s_0=0$. The general case follows similarly. Using Lemma \ref{neigh}
we can regard the equation (\ref{long}) as an equation in $(\ze,\eta) \in \Lol
\oplus \Lae $. More precisely, for a fixed $s$ small, we need to find $\ze \in
\Lol$ and  $\eta \in \Loe$ such that the approximate map  
$f= \exp _ {f_ \la } (   \ze+ P_\la g(  \ze) )$ solves the equation:
\be\label{ec0}
 \exp_f (P_ f \eta) = \exp_{f_\la} (s \xi) 
\end{equation}
The linearization of the equation (\ref{ec0}) at  $(0,\eta)$ is 
${\bold D} :  \Lol  \oplus \Lae \ra \Lo$, 
\[ {\bold D}_{(0,\eta)} ( \bz, \bn ) = \bz + P_\la \nabla g (\bz) +
 P_ \la  \bn + \Pi ( \bz,\eta )\]
\non where 
$\Pi ( \bz,\eta )$ is the derivative of $ P_\la \eta $ with respect to $ f_\la
$.
\bigskip

Our goal is to show that the operator $ {\bold D}_{(0,\eta)}$ is an isomorphism
in some appropriate norms on $ \Lol     \oplus \Lae $ and $ \Lo$.

\begin{defn}
 On $ \Lol     \oplus \Lae $ and $\Lo$ define the following norms:
\best
 \| (\bz,\bn) \|_{B_1}& =&\| \bz \lui + \| \bn +   \nabla g (\bz) 
    \lu  \;\;\mbox{ for any } \;\;(\bz,\bn) \in  \Lol\oplus\Lae\\
   \| \xi \|_{B_2}& =& \| D_\la \xi  \lu  \quad
 \mbox{for any} \quad \xi \in \Lo  
\eest
\end{defn}
 
\non Consider the operator 
$ {\bf T}: \Lol     \oplus \Lae \ra \Lo$ given by 
$ T( \bz, \bn) = z+ P_ \la ( \bn +   \nabla g (\bz) )$. 
Then ${\bf T}$ is continuous, since 
\best
 \| {\bf T}( \bz, \bn) \|_{B_2}& =& \| D_\la \bz + \bn + \nabla g (\bz) \lu 
\; \le \;  \| D_\la \bz  \lu  + \| \bn +   \nabla g (\bz)   \lu   \\ 
&\le&  
 C \la^ {1/4} \| \bz \lui + \| \bn + \nabla g (\bz) \lu\le \|(\bz,\bn)
\|_{B_1}
\eest

 \non for $ \la$  small enough. Recall that the low eigenvalues of 
$D_\la$ are of order $\la^{1/4}$, and thus 
$\| D_\la \bz  \lu  \; \le \;  \la^ {1/4}\| \bz \lui$ on $\Lol$.

\begin{lemma}\label{l10}
For $\la, \de $ small enough ${\bf T}$ is invertible, with the operator 
norm of the inverse uniformly bounded $\|  {\bf T}^{-1} \|  \; \le \; C_T$
(independent of $\la, \de $).
\end{lemma}
\pf Let $\alpha = \bz + P_\la( \bn +   \nabla g (\bz)   ) $. We need to
estimate $ \| \bz \lui $ and  $\| \bn +   \nabla g (\bz) \lp$ in terms of $ \|
\alpha  \|_{B_2}$. Since $D_ \la  \alpha =D_ \la  \bz + \bn+   \nabla g (\bz) $
then 
\best
 \| \bn+   \nabla g (\bz)   \lu   & \le&  
\| \alpha  \|_{B_2} + \|  D_ \la  \bz  \lu  \;   \le   \;
    \| \alpha  \|_{B_2} +C \la^ {1/4} \|  \ze  \lui   \\ 
&  \le & \| \alpha  \|_{B_2}  + C \la^ {1/4} \;\| \alpha  - 
P_\la (\bn+   \nabla g (\bz) )  \lui    \\ 
 & \le & \| \alpha  \|_{B_2} + C \la^ {1/4}\;  \| \alpha  \|_{B_2} +  
C \la^ {1/4}\;\| \bn+   \nabla g (\bz) \lu
\eest
 So for $\la$ small we get the uniform estimate 
$ \| \bn +  \nabla g (\bz)  \lu  \; \le \; C_1  \| \al \|_{B_2}$. This
gives  
\best 
 \| \bz  \lui & =&\| \alpha  - P_\la (\bn +   \nabla g (\bz) ) \lui  
\; \le \; 
 \| \alpha  \|_{B_2} + C \| \bn +   \nabla g (\bz) \lu \\
   &\le & C_2   \| \alpha \|_{B_2}
\eest
thus 
\best
 \| (\bz,\bn) \|_{B_1} \; \le \; C_T  \|  {\bf T}(\bz,\bn) \|_{B_2}
\eest 
So ${\bf T}$ is an injective linear operator. 
But by construction  ${\rm inded} ({\bf T})= 0$ thus ${\bf T}$ is invertible, 
with  $\|  {\bf T}^{-1} \|  \; \le \; C_T$ (independent of $\la, \de $).\qed

\begin{lemma}
For $\bz$ small, $ \| \Pi ( \bz, \eta) \|_{B_2}  \; \le \;  
C \| \eta  \lu   \| (\bz, 0 ) \|_{B_1}$.
\end{lemma}
\pf By differentiating the relation $D_f P_f \eta = \eta $ with respect to 
$f$ at $f_\la$ we get 
 \best \partial D_f ( P_ \la  \eta)(\bz) + D_ f  (\Pi ( \bz, \eta)) =0
\hskip.2in\mbox{ so}\\
\| \;D_ \la  (\Pi ( \bz, \eta)) \; \lu  = 
\| \;\partial D_f ( P_ \la  \eta)(\bz)\;  \lu   .
\eest
Using the expansion of 
\[ D_f \xi  = {1\over 2} \l( \nabla \xi + 
J(f) \circ \nabla \xi \circ j \r) + {1\over 8} N_f ( \partial_J f , \xi ) \]
(cf. \cite{MS}) it is easy to check that
\[  \|\; \partial D_f ( P_ \la  \eta)(\bz)  
 \; \lu  \; \le \; C \| \bz \|_{\infty,\la} \|  P_ \la  \eta   \lui  \]
 uniformly in a neighborhood of $ f_ \la $. Therefore 
\best
 \| \Pi ( \bz, \eta) \|_{B_2}   &=&  \|D_ f  (\Pi ( \bz, \eta))\; \lu  
 \; \le \; 
C \|\bz \|_{\infty,\la} \| P_ \la  \eta   \lui  \\
&\le&   C \| \bz  \lui   \| \eta  \lu  
 =  C \| \bz \|_{B_2}  \| \eta \lu.\hskip.3in \Box
\eest 

\non If we choose $\de$ small enough then for $ \| \eta  \lu  < \de$, 
\best
 \| \Pi ( \bz, \eta) \|_{B_2}   \; \le \; C_T/2 \; \|(\bz,\bn) \|_{B_2} 
\eest 
where $C_T$ is the constant in Lemma $\ref{l10}$ so 
$  {\bold D}_{(0,\eta)} ( \bz, \bn ) = {\bf T} ( \bz, \bn )+ 
\Pi ( \bz,\bn )  $ is still invertible. This concludes the proof of  
Proposition \ref{prop}. \qed

\subsection{The leading order term of the obstruction $\psi_t$ for $t$ small}
\label{Psi}

Next step is to identify the leading order term of the section 
$\psi_t$ as $t\ra0$. Let $\cn$  denote some stratum of $\ov\cu_d$ and 
$Gl\ra \cn$ denote the gluing data as  in (\ref{gl}). 
 For the sake of the gluing
construction, the gluing data has to be defined on the domain of the
bubble tree. But we will see in a moment that the important information
is encoded in the image curves. Introduce first some notation: 
 If  $u_i\in T_{y_i}S^2$ is a unit frame and $ \la_i$ is the gluing
parameter, let   
\[ v_i=\la_i\cdot u_i \in T_{y_i} S^2,\;\; (v_i\ne 0)\;\;\;\;
\mbox{ denote the gluing data.}\]

\begin{defn} For any $[f,y,v]\in Gl$, such that $f:B\ra \P^n$ is an element of
 $\cn$,
 let  $([f_i,y_i,v_i])_{i=1}^m$ be the bubble maps together
with the gluing data and let  $u$ be the image of the ghost base (so 
$u=f_j(y_j)$ for all $j\in\wt B$). Set
\bear
 \label{da}a([f,y,v])&=&
\ma(\sum)_{i=1}^n \ma(\sum)_{j\in \wt B}
\lg \;df_j(y_j)(v_j) \;,\; X_i\;\rg \om_i \\
\label{dbarn}  \bar \n(x)&= &
\ma(\sum)_{i=1}^n  \int_ {T^2} 
\lg \;\n(z, u) \;,\;\om_i(z)\;\rg  \om_i 
\eear
\non where $\{\om_i=X_idz,\; i=1,n\}$ is an orthonormal base of
$\ho(u)$, $X_i\in T_u\P^n$ and $\wt B$ is as in Definition \ref{obstrbd}. 
\end{defn}
 
Note that $a$ depends {\em only} on the gluing data on the first level $\wt B$
of
essential bubbles, and $\bar\nu$ depends only on the image of the ghost
base. Then
\begin{lemma}{\label{lg4}}
Using the notation above, let $\fl$ be an approximate gluing map. Then 
for $t$ and $|\la|=\sqrt{\la_1^2+\dots\la_i^2}$ small enough,
\bear\label{prn}
\pr ^{f_\la}(  \n )&=& \bar \n(u) +\o(|\la|)  \\    
\label{prfl}
 \pr ^{f_\la}( \del f_\la )&=& a([f,y,v]) + \o(|\la|^{3/2}). 
\eear    
and  the section $\psi_t$ has the form 
\bear\label{psi} 
\psi_t([f,y,v]) =t \bar\nu(u)+  a([f,y,v]) +
\o ( |\la| ^{5/4}+ t \sqrt{|\la|} + t^2). 
\eear
The estimates above are uniform on $\cn$.
\end{lemma}
\pf For the first 2 relations, it is enough to check them on components.
Assume for simplicity that $\wt B$ consists of
a single bubble $[f,y,v]$. If $\om=X dz$   
 is an element of the base for $H^{0,1}$, let $\omb$ be
the element of the local orthonormal frame for $\lal(f_\la) $ provided by
Lemma $\ref{lg2}$. Then 
\best
& |\lg  \n ,\oml-\omb \rgl |   \; \le \;  
 \|\n \|_\infty \| \oml-\omb \lu      \; \le \;   C \la& \quad \mbox{ so }\\
&\lg \n,\omb \rgl=  \lg \n,\oml \rgl+\o (\la)&
\eest
On the other hand, using the definition of $\oml$ and the fact that 
$\lg\;,\;\rgl$ is the usual inner product on $T^2$ off a small ball we get 
\best
 \lg \n,\oml \rgl&=& \ma(\int)_{|z|>\sla} \lg \n(z,f(y)), \om \rg =
\ma(\int)_ {T^2} \lg \n(z,f(y)), \om \rg + \o(\la)\;\;\mbox{ so}\\
 \lg \n,\omb \rgl&=& \ma(\int)_ {T^2} \lg \n(z,f(y)), \om \rg + \o(\la)
\eest
which gives (\ref{prn}). Similarly,
\[ |\lg  \del f_\la ,\oml-\omb \rgl |   \; \le \;   
\| \del f_\la \lu    \| \oml-\omb \lu      \; \le \;     C \la^{1/2}  \la 
  \; \le \;     C \la^{3/2} \]
\non and using the estimate (\ref{pe1}) and the definition of $\oml$ we get
\best \lg  \del f_\la ,\oml \rgl &=& 
\ma(\int)_{\sla   \; \le \;   |z|   \; \le \;   2\sla}  
{\sla \over |z|}\;  d \b  \lg  d f (y)(u), X \rg + \o(\la^2)\\ \\
&=&
\la \lg d f (y)(u), X \rg + \o(\la^2)   
\eest
Combine the previous 2 relations we get 
\best
 \lg  \del f_\la ,\omb \rgl = \lg d f (y)(\la u), X \rg + \o(\la^{3/2})
 \lg d f (y)(v), X\rg +\o(\la^{3/2})
\eest
 which implies  (\ref{prfl}).
\smallskip

The general case when $B$ has more bubbles follows in a similar maner
using the relation 
(\ref{aprglbd}) and the fact that $\om$ is 0 pass the first level of
nontrivial bubbles. 
\smallskip

Finally, the relation  (\ref{psi}) is a consequence of  (\ref{fla}) provided we
have an
estimate of the the quadratic part. For that use  (\ref{floer2}) to get 
\best
 \lg\; N_\la (P_\la \eta)\;,\; \omb \rgl   & \le &   
\| N_\la (P_\la \eta)\;  \|_{4/3,\la} \;\| \omb \|_{4,\la}   \; \le \;  
C \| \eta \|_{2,\la} \;\| \eta \|_{4/3,\la} \;\| \om \|_{4}   \\   
& \le &   \o( |\la|^{1/2} +t)\; \o (|\la|^{3/4}+t). 
\eest
\non Thus the  quadratic part is $\o( |\la| ^{5/4}+t\sqrt{|\la|}+ t^2)$. \qed
\bigskip

\step{ The definition of $\wt L\ra \wt\cu_d$ }
{}From this point on, since we are going to look at the leading order term, 
 it will become easier if we forget part of the gluing data.
We have already observed that the map $a$ depends only on the gluing
data on the first level $\wt B$ of essential bubbles. Moreover, if we 
denote by
\bear\label{w}
 w=\ma(\sum)_{j\in \wt B} df_j(y_j)(v_j)\in T_u\P^n
\eear
then the map $a$ and the linear part $\wt\psi_t$ of $\psi_t$ become 
respectively
\bear\label{ao}
a(w)&=&\ma(\sum)_{i=1}^n\lg w, X_i\rg \om_i\\
\label{psitdo}
 \wt\psi_t(w)&=&t \bar\nu(u)+  a(w)
\eear 
 Introduce a  space $\cw$ together with a projection $\pi:\cw\ra\ov \cu_d$
such that  
the fiber of $\pi$ at a 1-marked  curve (possibly with more components) 
is the span of
the tangent planes to all the image bubbles that meet at the marked
point. By definition $w\in \cw$ so (\ref{w}) defines a
projection  $p:Gl\ra \cw$.  Note though that $\pi:\cw\ra\ov\cu_d$ is not a 
vector bundle, and that $\cw$ it is equal to  the {\em relative tangent 
bundle} $L\ra \cu_d$ on the top strata of $\ov\cu_d$. 
\medskip

 Here is a more precise description of  $\cw$. Stratify $\ov\cu_d$ by 
letting $\cz_h$ be the union of all boundary strata
such that the image of the marked point is on $h$ nontrivial bubbles, i.e.
\bear\label{cz}
\cz_h=\{ f:B\ra\P^n \;|\; \wt B \mbox{ has } h \mbox{ elements }\}
\eear
Each $\ov\cz_h$ is a variety with normal crossings. 
For transversality arguments we need to use the 
moduli space $\czh_h$  obtained from $\cz_h$ by collapsing all the ghost 
bubbles up to the first level of essential bubbles. Note that 
 $\czt_{2}\supset\czt_3\supset\dots$, and the natural projection 
\[q:\cz_h\ra\czh_h\]
 has  fiber $\cu_{0,h}=\cm_{0,h+1}$, the moduli space of $h+1$ marked points on
the
sphere. Moreover,
\bear
{\rm dim}\; \cz_h= n-h-1\quad \mbox{ and } \quad
{\rm dim}\; \czh_h= n-2h+1
\eear 
In particular, $\cz_h\ne\emptyset$ only for  $h\le [{n+1\over 2}]$.  
\medskip

Let $L_i$ be the pullback of the relative tangent sheaf to the $i$'th 
factor of $\czh_h$. When the constraints  $\b_1,\dots,\b_k$ are in 
generic position, the fibers of  $L_1,\dots,L_h$ over a point in $\czh_h$
are linearly independent subspaces of $\P^n$. This is
because linear dependence imposes $n+1-h$ conditions, and 
$\czh_h$ is only $n-2h+1$ dimensional.  So on $\cz_h$ 
\bear\label{cw}
\cw|_{\cz_h}= q^*(L_1\oplus\dots\oplus  L_h)
\eear 
\rem Since not all the gluing parameters can be zero, a dimension
count argument similar to the one above shows that  $w$ defined by 
(\ref{w}) is an element of $\cw-\{0\}$, 
 the space nonzero vectors in $\cw$, thus  $p:Gl\ra \cw-\{0\}$.
\smallskip

Note that  $\cw|_{\cz_h}$ is nothing but the normal bundle of 
$\cz_h$ in $\ov\cu_d$, for any $2\le h\le [{n+1\over 2}]$. This
observation allows us to get a line bundle out of $\cw$ as follows:
\begin{defn}\label{ltilde}  Let $N=[{n+1\over 2}]$. 
Blow up $\ov\cu_d$ along $\cz_{N}$ (the bottom strata),
 then blow up the proper transform of  $\cz_{N-1}$ and so on, 
all the way up to blowing up the proper transform of $\cz_2$ and 
denote by 
\[\rho:\wt\cu_d\ra\cu_d\] the resulting
manifold. Similarly, after the first blow up, extend  $L$ over the 
exceptional divisor $E_{N}$ as the universal line bundle over 
$\P(N_{\cu_d}{\cz_N})$, the projectivization of the normal bundle of $\cz_N$,
 and so on. Let $\wt L\ra \wt \cu_d$ denote the blow up of $L$ constructed 
above. 
\end{defn}
By definition, the total space of $\wt L \ra \wt\cu_d$ is the
same as $\rho^*(\cw)$. From now on, we will make this identification.
\medskip

Both the map $a$ and the linear part $\wt\psi_t$ of $\psi_t$  pull back to 
$\wt L-\{0\}$ as
\bear\label{a}
a(w)&=&\ma(\sum)_{i=1}^n\lg w, X_i\rg X_i\\
\label{psitd}
 \wt\psi_t(w)&=&t \bar\nu(\pi(w))+  a(w)
\eear 
where $\pi:\wt L\ra \P^n$ is the composition $\wt L\ra
\wt\cu_d\ma(\longrightarrow)^{ev} \P^n$. For simplicity of notation, we 
have also  denoted by $\ev:\wt\cu_d\ra \P^n$ the composition 
$\wt\cu_d \ma(\ra)^\rho \cu_d\ma(\longrightarrow)^{ev}\P^n$.
 Note that by definition, $a$ is a
linear map but $\wt\psi_t$ is not, and we have the following diagramm:
\unitlength=1mm
\linethickness{0.4pt}
\begin{center}
\begin{picture}(95,15)
\put(20,12){\makebox(0,0)[cc]{$\wt L-\{0\}$ }}
\put(50,12){\makebox(0,0)[cc]{$\ev^*(T\P^n)$}}
\put(72,12){\makebox(0,0)[cc]{$T\P^n$}}
\put(35,0){\makebox(0,0)[cc]{$\cut_d$}}
\put(72,0){\makebox(0,0)[cc]{$\P^n$}}
\put(28,12){\vector(1,0){12}}
\put(34,15){\makebox(0,0)[cc]{$a,\;\wt \psi_t$ }}
\put(23,9){\vector(3,-2){9}}
\put(47,9){\vector(-3,-2){9}}
\put(70,9){\vector(0,-1){6}}
\put(72,3){\vector(0,1){6}}
\put(74,6){\makebox(0,0)[cc]{$\bar\nu$}}
\put(40,0){\vector(1,0){26}}
\put(55,2){\makebox(0,0)[cc]{$\ev$}}
\end{picture}
\end{center}
\begin{prop}\label{zeropsit}
As $t\ra 0$ the zero set of the section $\psi_t$ is homotopic to the
zero set of its leading order term 
\[\wt\psi_t:\wt L-\{0\}\ra \ev^*(T\P^n)\] 
\end{prop}
\pf
In generic conditions and for $t$ small enough the zero sets of both
sections 
\[\psi_t:Gl\ra\ev^*(T\P^n)\;\;\;\mbox{ and }\;\;\;
\rho^*p_*(\psi_t):\wt L-\{0\}\ra \ev^*(T\P^n)\]
consist of points lying on the top stratum of $\cu_d$ and $\cut_d$ 
respectively. But on the top stratum,  the projection 
$pr:Gl\ra\wt L-\{0\}$ is an isomorphism, thus the two 
 zero sets are diffeomorphic for $t$ small. Note that (\ref{psi}) gives
\[ p_*(\psi_t(w)) =t \bar\nu(u)+  a(w) + 
\o ( |w| ^{5/4}+ t \sqrt{|w|} + t^2)\] 
Finally,  Lemma \ref{nonvan} gives  that $w=\o(t)$ on the zero set of
$\psi_t$, so
\[  p_*(\psi_t(w)) =t \bar\nu(u)+  a(w) + \o ( t^{5/4})\]
giving the desired  homotopy as $t\ra 0$. \qed

\begin{lemma}\label{nonvan}
The linear map  $a:\wt L-\{0\}\ra \ev^*(T\P^n)$ defined in (\ref{a})
 has no zeros when the constraints 
$\b_1,\dots,\b_l$ are in a generic position, thus there exists $C>0$
such that
\bear\label{aw}
 |a(w)|\ge C|w|\eear
Moreover, there exists a uniform  constant $C$ on $\wt L-\{0\}$ such that the 
zero set of $\psi_t$ is contained in  $|w|\le Ct$.
\end{lemma}
\pf First part is a standard transversality argument and dimension
count. Note that $a$ induces a map
\best
a\otimes id:\wt L\otimes \wt L^*   &  \ra & \ev^*(T\P^n)\otimes \wt L^*  
\;\;\;\;\;\; \mbox{ i.e.}\\ 
a\otimes id:\ov \cu_d\ti \cx&\ra & \ev^*(T\P^n)\otimes \wt L^* 
\eest
Because of the $\cx^*$-equivariance of $a$, the zero set of 
$a:\wt L-\{0\}\ra \ev^*(T\P^n)$ is the same as the zero set of the section 
\best
 \wt a: \wt\cu_d&\ra&  ev^*(T\P^n)\otimes \wt L^* \\
 \wt a(x)&=&(a\otimes id) (x, 1)
\eest
If the constraints $\b_1,\dots,\b_k$ are in generic position, then 
 $\wt a$ is transverse to the zero set of $\ev^*(T\P^n)$. 
But the base $\wt \cu_d$ is only $n-1$ dimensional, while the fiber is $n$
dimensional, so generically $\wt a$ and thus $a$ has no zeros. 

 For the second part, note that on the zero set of $p_*(\psi_t)$ 
\best
0=p_*(\psi_t)&=& a(w)+t\bar\n(u)+\o(|w|^{5/4} +|w|^{1/2} \;t +t^2 )
\;\;\;\mbox{ so }\\
 a(w)&=&- t\bar\n(u)-\o(|w|^{5/4} +|w|^{1/2} \;t +t^2 )
\eest 
which combined with (\ref{aw}) gives
\best
C|w|\le  |a(w)|&\le& t|\bar\n(u)|+ \wt C(|w|^{5/4} +|w|^{1/2} \;t +t^2) 
\;\;\;\mbox{ i.e. }\\
|w|(C-\wt C|w|^{1/4})& \le &Ct
\eest
For $t$ and $w$ small, the left hand side is positive, completing the
proof. \qed

\subsection{The enumerative invariant $\tau_d$}\label{Tau}

Next step is to find 
the zero set of the leading order term of  $\psi_t$.
 As a warm-up we will discuss first the limit case $t=0$.

The constructions described in the previous sections 
 apply equally in this case, giving:
\begin{prop}\label{psit0}
Let $\cn$ be a ghost base boundary stratum of $\ov\cu_d$. 
Then the  moduli space of $J$-holomorphic tori 
close to $\cn$ is isomorphic to the zero set of a section in the
obstruction bundle over the space of gluing data
\[  \psi([f_i,y_i,v_i]_{i=1}^m) = a([f_i,y_i,v_i]]_{i=1}^m) 
+\o ( |\la| ^{5/4})\]
where a is defined by (\ref{da}).
  Moreover, for generic  constraints $\b_1,\dots,\b_l$, 
the number of $J$-holomorphic tori that define the enumerative invariant 
\[ \tau_d(\b_1,\dots,\b_l)\]
is finite, and the moduli space of these holomorphic tori is at a positive
distance
from the ghost base boundary strata of the bubble tree compactification. 
\end{prop}
\pf  For the second part, note that $\psi$ and $\la^{-1}\psi$ 
have the same zero set, so  
 as $\la\ra 0$ the limit of the end of the moduli space of $J$-holomorphic
tori is modeled by the zero set of the section $a$. But we have seen that
generically $a$ has no zeros, and thus there are no
$J$-holomorphic tori in a small neighborhood of that boundary stratum.\qed
\bigskip

Now we can now evaluate the contribution from the interior:
\begin{prop}\label{int}
For $t$ small, the number of $(J,t\nu)$-holomorphic maps that satisfy the 
constraints in the definition of $RT_{d,1}(\b_1\;|\;\b_2,\dots,\b_l)$
and are close to some $(J,0)$-holomorphic torus is equal to 
\[ n_j \tau_d(\b_1,\dots,\b_l)\]
where $n_j=|Aut_{x_1}(j)|$ is the order of the group of automorphisms of the
complex structure $j$ that fix the point $x_1$. 
\end{prop}
\pf Recall that $RT_{d,1}(\b_1\;|\;\b_2,\dots,\b_l)$ counts the number
of solutions of the equation
\[\del f(x)=\nu(x,f(x))\]
such that $f(x_1)\in \b_1$ and $f$ passes through $\b_2,\dots,\b_l$.
 
 A generic path of perturbations converging to 0 provides a cobordism
$\cm^\n$ to the solutions of the  equation 
\[\del f(x)=0\]
such that $f(x_1)\in \b_1$ and $f$ passes through $\b_2,\dots,\b_l$. 
A  $(J,0)$-holomorphic torus $f:T^2\ra \P^n$ is a smooth point of this
cobordism, i.e.  all the intersections are transversal and  
the cokernel $H^{0,1}(T^2, f^*(T\P^n))$ vanishes (since $f^*(T\P^n)$ is
a positive bundle for the standard complex structure). 
\smallskip

But the  invariant $\tau_d(\b_1,\dots,\b_l)$ counts the number of 
such solutions mod  the automorphism group of $j$. Imposing the
condition $f(x_1)\in \b_1$ reduces the stabilizer to just
$Aut_{x_1}(j)$.  \qed
\medskip 

\rem  Note that the pertubed invariant counts the number of $(J,\nu)$-
holomorphic maps with sign. This sign is  determined by the spectral flow 
of the linearization $D_f$ to $\del$. In the limit, when $\nu=0$, we have 
$D_f=\del$ thus all $(J,0)$-tori have a positive sign. This agrees with the 
way they were counted classically to obtain $\tau_d$. 
\medskip

\begin{lemma}\label{nonvan2}
For generic $\n$ the section 
$\bar\n:\ev_*(\ov\cu_d)\ra T\P^n $ defined by (\ref{dbarn})  has  
no zeros.
\end{lemma}
\pf  For generic $\n$, the section $\bar\n$ is transverse to the zero
section. But the fiber of $T\P^n$ is $n$ dimensional, and
the base $\im=\ev_*(\ov\cu_d)$ is only $n-1$ dimensional, so
$\bar\n$ has no zeros generically.\qed
\medskip

\rem The zeros $u\in\P^n$ of $\bar\n$  give the location of the point 
maps $u$ that can be perturbed away to get genus one  $(J,\n)$-holomorphic
maps representing $0\in H_2(\P^n)$. Since index=0 then
 generically $\bar \n$ has finitely many zeros. But $\im$ is a
codimension 1 subvariety in $\P^n$ that doesn't depend on $\n$. 
Then  we can
 choose  $\n$ generic so that its zeros do not lie in $\im$, and thus 
 $\bar \n(f(y)) \ne 0$ for any $[f,y]\in\ov\cu_d$.
\medskip

Moreover, Lemma \ref{zeropsit} showed that as  $t\ra 0$ the zero set $Z_t$
of $\psi_t$  is homotopic to the zero set $Z_0$ of the map
\best \psi_0: \wt L-\{0\}\ra \ev^*(T\P^n)\\ 
  \psi_0(w)=\bar \n(\pi(w)) +  a(w) 
\eest
where $a,\;\bar\n$ are defined in (\ref{a}), (\ref{dbarn}) and 
$\pi:\wt L\ra \P^n$ is the composition 
$\wt L\ra\wt\cu_d\ma(\ra)^{\ev} \P^n$. 
We have made a change of variables $w\ra w/t$.
\medskip

Next we identify the zero set $Z_0$.
Since $\bar \n(u)\ne 0$ on $\im$ then it induces a splitting of the 
obstruction bundle:
\bear\label{e} 
T\P^n\vert_{\im }=\cx<\bar \n> \oplus \;E 
\eear
\non where $E$ is an $n-1$ dimensional bundle, so 
\bear\label{splitev}
\ev^*(T\P^n) =\cx<\bar \n> \oplus\; \ev^* E 
\eear  
 \begin{lemma}\label{etil}
The number of zeros (counted with multiplicity)
 of $\psi_0$ is equal to
\[ c_{n-1}(\ev^*(E)\otimes \wt L^*) \]
\end{lemma}
\pf  Using  (\ref{splitev}) map $ \psi_0:\wt L-\{0\}\ra\ev^*(T\P^n) $ splits
as
\bear
\psi_1(w)=&\bar \n(\pi(w))\;\; +& a_1(w) \\
 \psi_2(w)=& &a_2(w)
 \eear
where $a_i$ denote the projections of $a(w)$. 
The map $a_2:\wt L-\{0\}\ra\ev^*(E)$ is $\cx^*$-equivariant, so 
tensored with the identity on $\wt L^*$ induces a $\cx^*$-equivariant map
\[ \bar a_2 :\wt\cu_d\ti\cx^*\ra \ev^*(E)\otimes \wt L^*\]
that has the same zero set as $a_2$. Let
\best
\wt a_2:\wt\cu_d\ra\ev^*(E)\otimes \wt L^* \;\;\;\;\mbox{ given by }\;\;\;
\wt a_2(x)=\bar a_2(x,1)
\eest 
Then  the zero set of $a_2$ is equal to $Z(\wt a_2)\ti \cx^*$. To find
the zero set of $\psi_0$, for any $(x,v)\in  Z(\wt a_2)\ti \cx^*$ solve
the equation
\[ 0=\psi_1(x,v)=\bar \n(x)+ a_1(x,v)=\bar \n(x)+v\cdot a_1(x,1)\]
Note that $a_1\ne 0$ on $Z(a_2)$ since $a$ has no zeros, so  for any 
$x\in Z(\wt a_2)$ there exists a unique $v\in \cx^*$ such that 
\[ -\bar \n(x)=v\cdot a_1(x,1)\]
This implies that there exists an isomorphism between the zero set of $\psi_0$
and
the zero set of $\wt a_2$. To complete the proof, note that for generic 
$\n$ the section $\wt a_2$ is 
transversal to the zero section of $\ev^*(E)\otimes \wt L^*$, so its zero
set is given by the Euler class of $\ev^*(E)\otimes \wt L^*$.\qed
\medskip

Finally, we can compute the boundary contribution:
\begin{prop}\label{zeros}
 For $t$ small, the number of $(J,t\nu)$-holomorphic maps that satisfy the 
constraints in the definition of $RT_{d,1}(\b_1\;|\;\b_2,\dots,\b_l)$
and are close to the boundary strata $\{x_1\}\ti\ov\cu_d$ is equal to 
\[  \ma(\sum)_{i=0}^{n-1}
{n+1\choose i+2}\ev^*(H^{n-1-i})\cdot c_1^i(\wt L^*)\] 
where $\wt L$ is the blow up of the relative tangent sheaf $L$ as 
in Definition \ref{ltilde}. 
\end{prop}
\pf As we have seen previously, the moduli space of 
$(J,t\nu)$-holomorphic maps that
 satisfy the constraints in the definition of
$RT_{d,1}(\b_1\;|\;\b_2,\dots,\b_l)$
and are close to the boundary strata $\{x_1\}\ti\ov\cu$ is diffeomorphic
to the zero set
of the section $\psi_0$. Using Lemma \ref{etil}, the later is equal to 
\[ c_{n-1}(\ev^*(E)\otimes \wt L^*) =\ma(\sum)_{i=0}^{n-1}
\ev^*(c_{n-i-1}(E)\;)\cdot c_1^i(\wt L^*)\]
But by definition $c_i(E)= c_{i}(T\P^n)={n+1\choose i}H^{i}$, completing
the  proof. \qed

\subsection{The other contribution}\label{Other}

 In the previous sections we have described in great length the gluing 
construction corresponding to the strata $\{x_1\}\ti \ov\cu_d$, that consists
of a ghost base and a bubble at the marked point $x_1$.  Finally, it is 
the time to sketch the gluing construction corresponding to other boundary
stratum $T^2\ti\ev^*(\b_1)$, and to explain why it does not give any 
contribution.
\begin{prop}\label{other}
 For $t$ small, the number of $(J,t\nu)$-holomorphic maps that satisfy the 
constraints in the definition of $RT_{d,1}(\b_1\;|\;\b_2,\dots,\b_l)$
and are close to the boundary strata $T^2\ti\ev^*(\b_1)$ is equal to 0. 
\end{prop}
\pf Construct first the space of approximate maps. The only difference
from the gluing construction decribed in Section \ref{Approx} is that we
need to allow the bubble point $x\in  T^2$ to vary. 
Since the tangent bundle of the torus is trivial, choose an isomorphism
\[ T T^2\cong T^2\ti\cx\]
which gives an identification $T _xT^2\cong \cx$ for all $x\in T^2$ 
(providing local coordinates on $T^2$). The set of gluing data will
then be modeled on:
\[ T^2\ti Fr\ti (0,\ep)\]  
where
\[Fr=\{\; [f,y,u] \;|\; [f,y]\in \ev^*(\b_1),\; u\in T_yS^2\;|u|=1\}\]
is the restriction of the frame bundle over $\cu_d$ defined by (\ref{fr}).
\medskip

To glue, use the unit frame $u\in TyS^2$ to identify $T_xT^2 \cong
T_yS^2$ which will induce natural coordinates on the sphere via the
stereographic projection. 
\medskip

Then all the constructions decribed in Sections \ref{Approx}-\ref{Tau} 
 extend to this
case. Since the holomorphic 1-form $\om\in H^{0,1}(T^2,\cx)$
is constant along the torus, then the isomorphism between the
obstruction bundle and $\ev^*(T\P^n)$ is independent of the bubble
point, so 
\[ \begin{array}{ccc}
\;\;\;H^{0,1}\cong\rho^*\ev^*(T\P^n)&&\ev^*(T\P^n)\\
\searrow\;\;\;\;\;\swarrow\;\;\;\;\;\;\;&& \downarrow\\
T^2\ti\ev^*(\b_1)\;\;\;&\ma(\longrightarrow)^\rho& \ev^*(\b_1)
\end{array}\]
Moreover, the linear part of the section $\psi_t$ that models the end of
the moduli space is also independent of the bubble point. 
But a dimension count shows that
the zero set of a $T^2$-equivariant section in the obstruction bundle
must be empty generically. \qed

\setcounter{subsection}{0} \setcounter{equation}{0}
\setcounter{theorem}{0} 
\section{Computations}\label{Comp}

In this second part of the paper we explain how one can compute the 
top power intersections $c_1^i(\wt L^*)\ev^*(H^{n-1-i})$ involved in 
Theorem \ref{gen}. The program is simple: first we  find recursive
formulas for the top intersections  $c_1^i( L^*)\ev^*(H^{n-1-i})$ 
 (Proposition \ref{pc1}), where
$L$ is the relative tangent sheaf of $\cu_d$, an object well known to the
algebraic geometers. Next we can  exploit the fact that $\wt L$ is a blow up
of $L$ to compute its corresponding top intersections.

Unfortunately, the
notation becomes quickly pretty complicated if we insist on keeping
track of all the information, so we chose to indicate at each step only the
new changes, leaving out the data that stays the same. 

\step{Notations}
If  $\b_0,\dots,\b_k$ are various codimension constraints let 
\[ \cu_d(\b_0\; ;\; \b_1,\dots,\b_k)=
\ev^*(\b_0)\; [\;\cu_d(\; ;\b_1,\dots,\b_k)\;]  \]
denote the  moduli space of 1-marked cuves in $\P^n$ passing through  
$\b_0,\dots,\b_k$, such that the special marked point is on
$\b_0$ and let
\[\cm_d(\b_0,\b_1,\dots,\b_k) \]
denote the corresponding  moduli space of curves (in which we 
forget the special marked point). 
\medskip

In particular, let $\cu_d= \cu_d(\; ;\b_1,\dots,\b_k)$ be the moduli 
space of 1-marked curves that appears in Theorem \ref{gen}. If  $i,j\ge 0$  
are  such that $i+j$= dim $\cu_d$ then  let 
\bear \phi_d(i,j\;|\;\b_1,\dots,\b_k)=
c_1^i(L^*)\;\ev^*(H^j)\;[\;\cu_d\;]
 \eear
denote the top intersection. Moreover, if $\wt \cu_d$ is the  
 blow-up  $\cu_d$ as in Definition \ref{ltilde}, let 
 \bear\label{xy}
 x=c_1(L^*) \in H^2(\cu_d,\Z),\;\;\wt x=c_1(\wt L^*) \in H^2(\wt\cu_d,\Z),
\;\; y=\ev^*(H)
\eear
where $y\in  H^2(\cu_d,\Z)$ or $y\in  H^2(\wt \cu_d,\Z)$ depending on the
context. Note that 
\bear
 \phi_d(i,j\;|\; \cdot \;)=x^i y^j\;[\;\cu_d\;]=x^i\;[\cu_d(H^j; \cdot )]
\eear

Using the notation above and the degeneration formula
(\ref{gendeg}), Theorem \ref{gen} becomes
\bear\label{bau}
n_j\tau_d(\cdot)=\ma(\sum)_{i_1+i_2=n}\si_d(H^{i_1},H^{i_2},\cdot)-
\ma(\sum)_{i=0}^{n-1}{n+1\choose i+2}\wt x^{i}y^{n-1-i}[\wt\cu_d]
\eear
\rem\label{steps} To compute a particular value for $\tau_d$ in $\P^n$ one
should use
a computer program based on the following four steps:
\begin{enumerate}
\item Find $\si_d$ using the recursive formula (\ref{spheres})
\item Find 
$\phi_d(i,j\; |\;  \cdot)=c_1^i(L^*)\ev^*(H^j)[\cu_d]$ 
using the recursive formulas of Proposition \ref{pc1}.
\item Find recursive formulas for 
$\wt x^i\cdot y^j=c_1^i(\wt L^*)\ev^*(H^j)[\wt \cu_d]$
as outlined in  Proposition \ref{pwtc1}.   
\item Finally, use (\ref{bau}) to get  $\tau_d$.
\end{enumerate}

\subsection{Recursive formulas for $c_1^i(L^*)\ev^*(H^{j}) $ }

Let $\cu_d$ be some $r$-dimensional moduli space of  1-marked curves of 
degree $d$ through some
constraints $\b_1,\dots,\b_k$ (not necessarily the same as in Theorem
\ref{gen}) and let $L\ra \cu_d$ be its relative
tangent sheaf.  In this section we give recursive
formulas for top intersections
\[ \phi_d(i,j\; |\;  \cdot )=c_1^i(L^*)\ev^*(H^j)[\cu_d] \]
where $i+j=r$ and 
 the constaints $\b_1,\dots,\b_k$ are dropped from the notation.

\begin{prop}\label{pc1} For every $r$-dimensional moduli space $\cu_d$ 
of any degree $d\ge 1$, there are the following recursive relations for the top
intersections: 
\bear\label{c1lini}
\phi_d(0,j\;|\; \cdot )&=&\si_d(H^j,\; \cdot ) \\
\nonumber
\phi_d(i+1,j\;|\; \cdot )&=&-{2\over d}\phi_d(i,j+1\;|\; \cdot \;)+
{1\over d^2}\phi_d(i,j\;|\;H^2,\; \cdot \;)\\ \nonumber 
&+&\ma(\sum)_{d_1+d_2=d\atop i_1+i_2=n}{d_2^2 \over d^2}\;
 \phi_{d_1}(i,j\;|H^{i_1}\;,\; \cdot \;)\cdot\si_{d_2}(H^{i_2},\; \cdot
\;)\nonumber\\ \label{c1li}
&+&\ma(\sum)_{d_1+d_2=d\atop i_1+i_2=n+j}{d_2^2 \over d^2}
 \phi_{d_1}(i-1,i_1\;|\; \cdot \;)\cdot\si_{d_2}(H^{i_2},\; \cdot \;)
\eear
for any $i\ge 0$, where the sums above are  over all possible distributions of

the constraints $\b_1,\dots,\b_k$ on the two factors and $d_1,d_2\ne 0$. 
When $i=0$, the last term in (\ref{c1li}) is missing. 
\end{prop}
\pf The first relation follows by definition, and  provides the initial 
step of the recursion. The second one requires more work. 
In what follows, we will identify a cohomology class like $c_1(L)$ with a
divisor representing it. Then:
\begin{lemma}\label{lc1l}
On $\cu_d$, we have the following relation:
\bear\label{c1l}
 c_1(L^*)={1\over d^2} {\cal H}-{2\over d}\ev^*(H) +{1\over d^2}
\ma(\sum)_{d_1+d_2=d}d_2^2 \cm_{d_1,d_2}
\eear
where ${\cal H}$ denotes the extra condition that the curve passes through
$H^2$, and  $ \cm_{d_1, d_2}$ 
denotes the  boundary stratum corresponding to the
splittings in a degree $d_1$ 1-marked  curve and 
 a degree $d_2$  curve, for $d_i\ne 0$ (for  all possible distributions
of the constraints $\b_1,\dots,\b_k$ on the two components).   
\end{lemma}
\pf Fix 2 hyperplanes in generic position in $\P^n$. 
Each curve in $\cu_d$ intersects a hyperplane in $d$ points. Then the
moduli space $Y=\ev^*_{k+1}(H)\cap\ev^*_{k+2}(H)$ 
of 1-marked  curves passing through $\b_1,\dots,\b_k,H,H$ is a  $d^2$ fold 
cover of $\cu_d$:
\[ \pi:Y \ra \cu_d, \hskip.3in 
[f,y_1,\dots,y_k,a,b\;;y]\ra [f,y_1,\dots,y_k\;;y]\]
 Define the section 
\[ s([f,y_1,\cdots,y_k,a,b\;;y])={(a-b)dy\over (y-a)(y-b)}\]
Then $s$ is a section in the relative cotangent sheaf $L^*$, and it
extends to the compactification $\ov{\cu_d}$. As $a$ and $b$ are getting
closer together, the section $s$ converges to 0. Thus its  zero set is the
sum of the divisors $\{a=b\}$ and $\cm(y\; ; \;a,b)$, where 
$\cm(y\; ;\;a,b)$ is the sum of all boundary strata  corresponding to
splittings into a degree $d_1$ 1-marked  bubble and a degree $d_2$ bubble
containing $a,b$ for $d=d_1+d_2$. Note that $d_i\ne 0$. 
The infinity divisor is $\{y=a\}+\{y=b\}$.
Thus
\[\pi^*( c_1(L^*))=\; \{a=b\}+ \cm(y\; ; \;a,b)-\{y=a\}-\{y=b\}\]
Note that  
\[ d^2 c_1(L^*)=\pi_*\pi^*( c_1(L^*))\]
When projecting  down to $\ov \cu_d$, the divisor $\{a=b\}$ becomes 
 ${\cal H}$, and the divisors   $ \{y=a\}$, $ \{y=b\}$ become each 
$d\cdot\ev^*(H)$. The rest amounts to summing over all codimension 1 boundary
strata. The boundary strata $\cm_{d_1, d_2}$
appears with coefficient $d_2^2$ in $\pi_*(\cm(y,a\; ; \;b))$. 
Combining all the pieces  together completes the proof of Lemma. \qed
\medskip

\rem We could have chosen any 2 marked
points out of the already existent ones, and then express $c_1(L)$ in terms
of them. But then this  expression would not look independent of choice. 
Nevertheless, with some work, one can actually see that all these divisors
are homotopic. We have chosen to introduce 2 new marked points to avoid
this issue.
\medskip

\rem Note that the base locus of the line bundle $L$ is exactly the
union of the divisors $\cz_h$.  Doing the intersection theory in the
blow up along the base locus is the same as considering the excess
intersection (see \cite{ful}).
\bigskip

Relation  (\ref{c1l}) provides the basic relation for proving
(\ref{c1li}): 
\[ c_1^{i+1}(L^*)=-{2\over d}c_1^i(L^*)\cdot \ev^*(H) + 
{1\over d^2}c_1^i(L^*)\cdot  {\cal H}+
\ma(\sum)_{d_1+d_2=d}{d_2^2 \over d^2}\;
c_1^i(L^*)\cdot \cm_{d_1,d_2} \]
so taking a cup product with $\ev^*(H^j)$ we get:
\bear\label{c1ini}
\phi_d(i+1,j\;|\; \cdot )&=&-{2\over d}\phi_d(i,j+1\;|\; \cdot )+
{1\over d^2}\phi_d(i,j\;|\;H^2,\; \cdot )\\ \nonumber
&+&\ma(\sum)_{d_1+d_2=d}{d_2^2 \over d^2}\;
c_1^i(L^*)\ev^*(H^j)\cdot \cm_{d_1,d_2} 
\eear
 Next, we need to understand the restriction of $L^*$ to the
boundary stratum $\cm_{d_1,d_2}$.
Let 
\[p:\cm_{d_1,d_2}\ra \cu_{d_1}\] 
be the projection on the first
component (the one that contains the special marked point $y$). If $A,B$ are 
 the 2 special points of $\cm_{d_1, d_2}$ (where the 2
components meet), let
\[\ev_A\ti\ev_B:\cu_{d_1}\ti \cm_{d_2}\ra \P^n\ti \P^n\]
be the corresponding evaluation map. Then by definition
\bear\label{cm12}
\cm_{d_1,d_2}=(\ev_A\ti\ev_B )^*([\De])
\eear
 where $\De$ is the diagonal of $\P^n\ti\P^n$. 
Moreover, it is known that as divisors,
\bear\label{c1la}
c_1(L^*)/\cm_{d_1,d_2} =p^* c_1(L^*_A)+\{ y=A\}
\eear
where $L_A=L_{|\cu_{d_1}}$ is the relative tangent sheaf of
$\cu_{d_1}$. Next step is to find 
\bear\label{c1cmd1d2} 
 c_1^i(L^*)/\cm_{d_1,d_2}=\ma(\sum)_{l=0}^i {i\choose l}
\;p^* c_1^{i-l}(L^*_A)\cdot(\{ y=A\})^l  
\eear 
For the  self intersection of the divisor $\{ y=A\}$ note that its
normal bundle $N$ inside $\cm_{d_1,d_2}$ is nothing but 
$p^*(L_A)/\{ y=A\}$, so for $l>0$,
\[ (\{ y=A\})^l= c_1(N)^{l-1}=(-1)^{l-1}p^*c_1^{l-1}(L_A^*)
\cdot [\{y=A\}]  \]  
Substituting in  (\ref{c1cmd1d2}) and after some algebraic manipulations
we get:
\bear\label{pof}
c_1^i(L^*)\cdot[\cm_{d_1,d_2}] =p^* c_1^i(L^*_A)+p^*c_1^{i-1}(L_A^*)\cdot
[\{ y=A\}]
\eear

We will do the intersection theory inside $\cu_{d_1}\ti\cm_{d_2}$. The
relation  (\ref{cm12}) combined with  
 $[\De]=\ma(\sum)_{i_1+i_2=n}H^{i_1}\ti H^{i_2}$  gives
\best
\ev^*(H^j)\cdot[\cm_{d_1,d_2}]&=&\ma(\sum)_{i_1+i_2=n}
\cu_{d_1}(H^j\; ;\;H^{i_1},\; \cdot \;)\ti \cm_{d_2}(H^{i_2},\; \cdot \;)\\
\ev^*(H^j)\cdot[\{ y=A\}]&=&\ma(\sum)_{i_1+i_2=n+j}
\cu_{d_1}(H^{i_1}\;;\; \cdot \;)\ti \cm_{d_2}(H^{i_2},\; \cdot \;)
\eest
where we sum over all possible distributions of the constraints on the
two components. The relations above imply
\bear
\nonumber 
&&\hskip-.6in \ev^*(H^j)\cdot p^*c_1^i(L^*_A)\cdot [\cm_{d_1,d_2}]\\ \nonumber
&&=\ma(\sum)_{i_1+i_2=n}(c_1^i(L_A)\cdot 
\cu_{d_1}(H^j\; ;\;H^{i_1},\; \cdot \;)\;)\ti
 \cm_{d_2}(H^{i_2},\; \cdot \;)\\ 
&&=\ma(\sum)_{i_1+i_2=n}
\phi_{d_1}(i,j\; |\;H^{i_1},\; \cdot \;)\cdot
 \si_{d_2}(H^{i_2},\; \cdot \;)
\eear 
\bear\nonumber
&&\hskip-.6in\ev^*(H^j)\cdot p^*c_1^{i-1}(L_A^*)\cdot[\{ y=y_A\}]
\\&&=\ma(\sum)_{i_1+i_2=n+j} \nonumber
(c_1^{i-1}(L_A)\cdot\cu_{d_1}(H^{i_1}\; ;\; \cdot \;)\;)\ti
 \cm_{d_2}(H^{i_2},\; \cdot \;)\\
&&=\ma(\sum)_{i_1+i_2=n+j}
\phi_{d_1}(i-1,i_1\; |\; \cdot \;)\cdot \si_{d_2}(H^{i_2},\; \cdot \;) 
\eear 
Substituting these relations in (\ref{c1ini}) using  (\ref{pof}) we get
(\ref{c1li}), which concludes the proof of Proposition \ref{pc1}.  \qed

\subsection{Recursive formulas for  $c_1^i(\wt L^*)\cdot\ev^*(H^{j})$ }

Next step is to express the top intersections involving the first Chern
class of $\wt L$, the blow up of $L$, in terms of the top intersections
involving the first Chern class of $L$. The program for such kind of
computations is very nicely outlined in \cite{ful}, which we will follow
closely. Although recursive formulas can be found for any $n$, the more
strata we need to blow up, the longer and more complicated looking 
these formulas become. 

For simplicity of the presentation, in this section we will give only
the general principles of the algorithm, without working out completely
the recursive formulas. In the next section we will use this
algorithm to obtain recursive formulas for small values of $n$ (i.e.
$n\le 4$).  
\medskip

Let $\cu_d=\cu_d(\; ;\b_1,\dots,\b_k)$ the some $r$-dimensional moduli 
space of 1-marked  curves.  Recall the construction of  $\cut_d$: 
starting with $\cu_d$, we first blow up  along $\cz_{N}$, then we blow up 
the proper
transform of $\cz_{N-1}$ and so on, up to blowing up the
proper transform  of $\cz_2$. Since $\wt L$ extends as the blow up of $L$
then 
\bear
 c_1(\wt L)= c_1(L)+\ma(\sum)_{h=2}^{N} E_h
\eear
where $E_h$ is the exceptional divisor corresponding to the proper
transform of $\cz_h$.

\begin{prop}\label{pwtc1} Using the notations above, the top intersections
\bear
c_1^i(\wt L^*)\ev^*(H^{j})
\eear
on $\wt\cu_d$ can be recursively expressed in terms of the top intersections 
\best
\phi_d(k,l\;|\;\cdot\;)=c_1^k(L^*)\cdot\ev^*(H^{l})
\eest on possibly lower dimensional moduli spaces $\cu_{d'}$. 
\end{prop}
\pf The idea is of course to do inductively one blow up at a time. 
Although the fact there exist such recursive formulas is not that hard
to see, writting them down becomes pretty complicated very quickly. So
we explain why such formulas exist, leaving their  derivation 
for later. By definition,
\best
\wt x=c_1(\wt L^*)=x-\ma(\sum)_{h=2}^{N} E_h
\eest
Let
\bear
x(h)=x-\ma(\sum)_{l=h}^N E_l
\eear
so $\wt x(N+1)=x$ and $\wt x(2)=\wt x$. Using $x(h)=x(h+1)-E_h$, and  
expanding, 
\bear\nonumber
x(h)^i\cdot y^{j}&=&x(h+1)^{i}\cdot y^{j}+
\ma(\sum)_{l=1}^i{i\choose l}x(h+1)^{i-l} y^{j}(-1)^l\; E_h^l \\
&=&x(h+1)^{i} y^{j}-
\ma(\sum)_{l=h}^i{i\choose l}x(h+1)^{i-l} y^{j} 
s_{l-h}(N_{\wt \cu_d}\czt_h) [\czt_h]
\label{wtxiyj}
\eear
The last equality is a consequence of the following 
\begin{fact}(cf. \cite{ful})  Assume  $X\subset Y$ regular 
imbedding of codimension $a$, where dim $Y$=r. Let $\pi:\wt Y\ra Y$ be the
blow up of $Y$ along $X$ and $E=\P(N_YX)$ be the exceptional curve. For any 
$\al\in H^{r-l}(Y)$, the top intersection 
\bear
(\pi^*(\al)\cup E^l)\cap[\wt Y]=(-1)^{l-1} (\al\cup s_{l-a}(N_YX))\cap[Y]
\eear
as integers, where $l\ge 1$ and $s(N)$ is the segree class of the normal
bundle $N$. 
\end{fact} 

Next we need to understand how $x(h+1)$ or equivalently $E_m$ restricts
to $\czt_h$, and also we need to find $s(N_{\wt\cu_d}\czt_d)$. First, we
find:
\step{The normal bundle of $\cz_h$ in $\cu_d$} 
Recall that  $\cz_h$ consists of bubble trees with $h$ essential
components meeting at the image of the ghost base. So in particular, 
$\cz_h$ has components indexed by the different distributions of the degree
on the  $h$ bubbles:
\bear\label{ovcz}
\cz_{d_1,\dots,d_h}=\ev_0^*([\De])\subset 
\ov\cm_{0,h+1}\ti  \cu_{d_1}\ti\dots\ti \cu_{d_h}
\eear
where $d_i\ne 0$ for $i=1,\dots,h$, $\De$ is the small 
diagonal in $(\P^n)^h$ and 
\bear
&&\ev_0:  \cu_{d_1}\ti\dots\ti \cu_{d_h}\ra (\P^n)^h\\
&&ev_0([f_1,y_1],\dots,[f_h,y_h])=(f_1(y_1),\dots,f_h(y_h))
\eear
is the evaluation map. Then
\bear
 \cz_h={1\over h!}\ma(\bigcup)_{d_1+\dots+d_h=d}\cz_{d_1,\dots,d_h}
\eear
where the factor of $h!$ comes from the action of the symmetric group 
that permutes the order of the $h$ bubbles (yielding the same bubble
tree).  Let  
\best
\ov \cm_{0,h+1}\ti \cu_{d_1}\ti\dots\ti  \cu_{d_h}&
\ma(\longrightarrow)^{p_i} &\cu_{d_i}
\eest
be the projection and $L_i$ be the relative tangent sheaf  of 
the $i$'th factor. It is easy to check that 
\begin{lemma}\label{normalbd}
Using the notations above, the normal bundle $N_{ \cz_h} \cu_d$ of 
$\cz_h$ in $\cu_d$ is isomorphic to 
\bear\label{normbd}
 p_1^*L_1\oplus \dots\oplus p_h^* L_h \quad \mbox{on each component}
\eear
so 
\bear
s(N_{ \cz_h} \cu_d)={1\over (1-x_1)\cdot \dots\cdot (1-x_h)}
\eear
where $x_i=p_i^*c_1(L_i^*)$. 
\end{lemma}
\medskip

\rem One word of caution: so far we have defined $x_i$ on each component
$\cz_{d_1,\dots,d_h}$  of $\cz_h$, but these definitions {\em do not 
match} on the intersection of two components. Nevertheless, after we blow
up $\cu_d$ as in Definiton \ref{ltilde}, all the components of $\cz_h$ 
become disjoint, so doing the intersection theory in the blow up allows
up to  treat each component separately as if they were disjoint.  
\medskip

Next, use the following 
\begin{fact}(cf. \cite{ful}) Assume  $X,\;Y\subset Z$ are regular 
imbeddings. Let $\wt Z=Bl_X Z$ be the blow up of $Z$ along $X$ and 
$\wt Y=Bl_{X\cap Y} Y$ be the proper transform of $Y$. 
Denote by $E=\P(N_Z X)$ the exceptional curve in $\wt Z$, and let 
$F=\P(N_Y {(X\cap Y})$ be the exceptional curve in $\wt Y$. Then:
\bear
E\cap \wt Y&=&F\\
N_{\wt Z} \wt Y&\cong& \pi^*(N_Z Y) \otimes \co(-F)\quad \mbox{ so } \\
\label{spN}
s_p(N_{\wt Z} \wt Y)&=&\ma(\sum)_{i=0}^{p-a} {a+p\choose a+i}
\;s_i(N_Z Y)\; F^{p-i}
\eear 
where $a=$rank $N_Z Y$.
\end{fact}
So
\bear
E_l\cap \czt_h=E_{h,l}
\eear
is the exceptional curve in the blow up of
$\cz_h$ along the proper transform of $\cz_h\cap\cz_l$ and 
\bear
N_{\wt \cu_d} \czt_h&\cong& \pi^*(N_{\cu_d} \cz_h) \otimes 
\co\l(-\ma(\sum)_{l=h+1}^N E_{h,l}\r)
\eear
Note that $\cz_h$ has less dimensions than $\cu_d$ and it 
is stratified by subvarieties $\cz_h\cap\cz_l$ for $l\ge h+1$ 
the same way $\cu_d$ is stratified by the sets $\cz_l$. Thus
we can repeat the same construction.
\medskip

Inductively, the intersection theory takes place inside a 
strata of form $\cz_{h_1}\cap\dots\cap \cz_{h_l}$. Now, it is not that
easy to list and parametrize all the possible bubble tree configurations
in this intersection. But a closer look reveals that although the 
combinatorics involved is complicated, all these intersections  have 
components of the form
\bear\label{piece}
ev_0^*([\De])\subset\cz\ti \cu_{d_1}\ti \dots\ti \cu_{d_m}
\eear
where $\cz$ is some substrata of $\ov \cm_{0,m+1}$, and $\De$ is the
small diagonal in $(\P^n)^m$. Such component comes in with a coeficient of one
over the order of the subgroup of permutations that preserve the same 
bubble tree configuration. 

The  normal bundle of (\ref{piece}) has the same form as in 
Lemma \ref{normalbd}, and using the arguments outlined above 
we are inductively  decreasing either the
number of exceptional curves or the dimension of the moduli space we do
the intersection theory over. In either case, the process terminates in
finite time, reducing the top intersections $\wt x^i y^j$ on $\wt\cu_d$ 
we started with to sums of top intersections of the form
\bear\label{mess}
x_1^{i_1}\dots x_m^{i_m} y^j \quad \mbox{ on }\quad
\ev_0(\De)\subset \cz\ti \cu_{d_1}\ti \dots\ti \cu_{d_m}
\eear
Finally, let 
\best
\pi:\cz\ti \cu_{d_1}\ti \dots\ti \cu_{d_m}\ra
\cu_{d_1}\ti \dots\ti \cu_{d_m}
\eest
be the projection. Since all the classes in (\ref{mess}) are pull-backs
by $\pi$ then the top intersection (\ref{mess}) vanishes unless
$\cz\subset \cm_{0,m+1}$ is 0 dimensional. When $\cz$ is 0-dimensional
then using the decomposition of the diagonal 
\best
[\De]=\ma(\sum)_{i_1+\dots+i_m=n}H^{i_1}\ti\dots\ti H^{i_m}
\eest
and letting $y_j=\ev^*(H)$ on the $j$'th factor we get
\bear
 \ev_0^*([\De])&=&\ma(\sum)_{j_1+\dots+j_m=n}
y_1^{j_1}\dots y_m^{j_m}
[\cu_{d_1}\ti\dots\ti\cu_{d_m}]\quad \mbox{ and}\\
y^j \ev_0^*([\De])&=&\ma(\sum)_{j_1+\dots+j_m=n+j}
y_1^{j_1}\dots y_m^{j_m}[\cu_{d_1}\ti\dots\ti \cu_{d_m}]
\eear
so (\ref{mess}) is equal to
\bear
\ma(\sum)_{j_1+\dots+j_m=n+j}
(\;x_1^{i_1}y_1^{j_1}[\cu_{d_1}]\;)\ti\dots \ti
(\;x_m^{i_m}y_m^{j_m}[\cu_{d_m}]\;)\\ \label{splity}
=\ma(\sum)_{i_1+\dots+i_m=n+j}
\phi_{d_1}(i_1,j_1)\cdot\dots\cdot \phi_{d_m}(i_m,j_m)
\eear
giving Proposition \ref{pwtc1}. \qed

\section{Applications to $\P^n, \;n\le 4$}

Finally, we apply the inductive algorithm described in the
previous section to obtain recursive formulas for the elliptic 
 enumerative invariant $\tau_d$ in $\P^n$, for $n\le 4$. In this case
the story is quite simple, since  $N=\left[ {n+1\over 2}\right]\le 2$, 
so we need to blow up at most one strata, $\cz_2$.
 
\step{Explicite formulas for  $c_1^i(\wt L^*)\cdot\ev^*(H^{j})$}
 Note that for $n=2$ there is nothing to 
blow up, so
\bear\label{ln2}
\wt L=L  \quad \mbox{ for }\quad n=2
\eear  
A litlle more work gives:
\begin{lemma}\label{l34} Using the notations in  Theorem \ref{gen}, when 
$n=3,\;4$, we have the following relations:
\bear\label{x1y}
&&\hskip-.3in c_1(\wt L)\cdot \ev^*(H^{n-2})=\phi_d(1,n-2\;|\;\cdot\;)\\
\label{x2y}
&&\hskip-.3in c_1^2(\wt L)\cdot \ev^*(H^{n-3})=\phi_d(2,n-3\;|\;\cdot\;)-
\hskip-.2in\ma(\sum)_{d_1+d_2=d\atop i_1+i_2=2n-3}\hskip-.1in
\si_{d_1}(H^{i_1})\cdot  \si_{d_2}(H^{i_2})\\
\label{x3y}
&&\hskip-.3in c_1^3(\wt L)=\phi_d(3,0\;|\;\cdot\;)-
\ma(\sum)_{d_1+d_2=d\atop i_1+i_2=4}
\phi_{d_1}(1,i_1\;|\;\cdot\;)\cdot \si_{d_2}(H^{i_2})
 \mbox{ for }n=4
\eear
where the sums above are over all possible distributions of the
constraints on the two components and $d_i\ne 0$. 
\end{lemma}
\pf When $n=3,\; 4$ we need to blow up only $\cz_2$. Use (\ref{wtxiyj}) 
to get:
\bear\nonumber
\wt x^i y^{n-1-i}\;[\wt\cu_d]&=& x^{i}y^{n-1-i}\;[\cu_d] \\ \label{blxy}
&-&\ma(\sum)_{l=2}^i
{i\choose l}x^{i-l}s_{l-2}(N_{\cz_2}) y^{n-1-i} \; [\cz_2] 
\eear
Note that for $i=1$ the sum in (\ref{blxy}) is indexed by the empty set, thus
giving (\ref{x1y}). When $i=2$ the sum reduces to:
\best
x^{0}\cdot s_{0}(N_{\cz_2})\cdot y^{n-3} \; [\cz_2]= y^{n-3} \; [\cz_2]=
{1\over 2}\ma(\sum)_{d_1+d_2=d\atop i_1+i_2=2n-3} \si_{d_1}(H^{i_1})
\cdot \si_{d_2}(H^{i_2})
\eest
by (\ref{splity}), giving  (\ref{x2y}).

When $n=4$ and $i=3$ then the sum in (\ref{blxy}) becomes
\bear\label{43}
(3 x\cdot s_0(N_{\cu_d}{\cz_2})+s_1(N_{\cu_d}{\cz_2}))[\cz_2]
\eear
But note that $x|_{[\cz_2]}=0$ and Lemma \ref{normalbd} gives
\best
s(N_{\cu_d}{\cz_2})={1\over (1-x_1)(1-x_2)}\quad \mbox{ thus } \quad 
s_1(N_{\cu_d}{\cz_2}))=x_1+x_2
\eest
So (\ref{43}) becomes
\best
(x_1+x_2)\;[\cz_2]&=&
{1\over 2}\ma(\sum)_{d_1+d_2=d\atop j_1+j_2=n}(x_1+x_2)y_1^{j_1}y_2^{j_2}
\;[ \cu_{d_1}\ti\cu_{d_2}]\\
&=&\ma(\sum)_{d_1+d_2=d\atop j_1+j_2=n}(x_1y_1^{j_1} [\cu_{d_1}])\cdot
(y_2^{j_2} [\cu_{d_2}])= \ma(\sum)_{d_1+d_2=d\atop j_1+j_2=n}
 \phi_{d_1}(1,j_1\;|\;\cdot\; )\cdot\si_{d_2}(H^{j_2},\;\cdot\;)
\eest
using again (\ref{splity}). \qed
\medskip

Now we can prove for example that:
\begin{prop}\label{p^2}
The number $\tau_d(p^{3d-1})$  of degree $d$ elliptic 
curves  in $\P^2$ with 
fixed $j$ invariant and passing though $3d-1$ points is 
\bear 
\tau_d(p^{3d-1})={2\over n_j}{d-1\choose 2}\si_d(p^{3d-1})
\eear
where $\si_d$ is the number of rational curves through $3d$ points,
 and $n_j$ is the order of the group of automorphisms of the complex structure
$j$ fixing a point. 
\end{prop}
\pf For $n=2$, relation (\ref{conden}) combined with (\ref{ln2}) gives: 
\bear\label{relp2}
n_j \tau_d(p^{3d-1})= \si_d( l,l,p^{3d-1})-3\ev^*(H)-c_1( L^*)
\eear
where $L\ra \cu_d$ is the relative tangent sheaf over the moduli
space of 1-marked  rational curves of degree $d$ passing through $3d-1$. 
 The moduli space $\cm_d$ of unmarked curves is $n-2=0$ dimensional, 
consisting of $\si_d(p^{3d-1})$ curves. Using (\ref{pc1}) 
(or easier by inspection)
 \[ c_1(L^*)=-{2\over d}\si_d(l,p^{3d-1})=-2\si_d(p^{3d-1})\]
\[\ev^*(H)=\si_d(l,p^{3d-1})=d\si_d(p^{3d-1})\;\;\;\mbox{ and }\;\;\;
 \si_d(l,l,p^{3d-1})=d^2 \si_d(p^{3d-1})\]
So plugging them back in (\ref{relp2}) we obtain
\[ \tau_d(p^{3d-1})={1\over n_j}(d^2-3d+2)\;\si_d\]
which gives (\ref{p^2}). \qed
\medskip

In particular, 
\vskip-.5in
\bear 
\tau_d(p^{3d-1})= \left\{\begin{array}{ll}
{d-1\choose 2}\si_d& \mbox{ if } j\ne 0, 1728\\ \\
{1\over 2} {d-1\choose 2}\si_d& \mbox{ if } j= 0\\ \\
{1\over 3} {d-1\choose 2}\si_d& \mbox{ if } j=1728
\end{array}\right.
\eear
This formula was recently obtained  by Panharipande  \cite{rp} using
 different methods.
\medskip

\non Next we can prove that:
\begin{prop}\label{compp3} The number $\tau_d=\tau_d(p^a,l^b)$ of elliptic
curves 
 in
$\P^3$ with fixed $j$ invariant and passing through $a$ points and $b$
lines (such that $2a+b=4d-1$) is given by:
\bear\label{formp3}
\tau_d={2(d-1)(d-2)\over dn_j}\si_d(l)-{2\over dn_j}\ma(\sum)_{d_1+d_2=d}
d_2(2d_1d_2-d)\si_{d_1}(l)\si_{d_2}
\eear
where $\si_d(l)=\si_d(p^a,l^b,\;l)$ is the number of degree d rational 
curves in $\P^3$
passing through same conditions as $\tau_d$ plus  one more line. By the   
term $\si_{d_1}(l)\si_{d_2}$ we understand the sum over all
decompositions into a degree $d_1$ and a degree $d_2$ bubble such that the
constraints are distributed in all possible ways on the bubbles, and
$d_i\ne 0$.
\end{prop}
\pf When $n=3$, Theorem \ref{gen} gives:
\bear\nonumber
n_j\tau_d(p^a, l^b) &=&\ma(\sum)_{i_1+i_2=3} 
\si_d(H^{i_1},H^{i_2},p^a,l^b )- 6\ev^*(H^2)\\  \label{3ci}
&-&4\ev^*(H) c_1(\wt L^*)-c_1^2(\wt L^*)
 \eear 
The moduli space $\cm_d$ of degree $d$ unmarked curves  passing through $a$
points and $b$ lines is $n-2=1$ dimensional, with a finite
number of bubble trees in the boundary. Then Proposition \ref{compp3} is
a consequence of (\ref{3ci}) and the following
\begin{lemma} In $\P^3$,
\best
 &&\ma(\sum)_{i_1+i_2=3} \si_d(H^{i_1},H^{i_2},p^a,l^b )=2d
\cdot\si_d(p^a,l^{b+1}) \\ 
&&\ev^*(H^2)=\si_d(p^a,l^{b+1}) \\
 &&\ev^*(H)\cdot c_1(\wt
L^*)=\ev^*(H)\cdot c_1( L^*)=
 -{1\over d}\si_d(l)+
{1\over d}\ma(\sum)_{d_1+d_2=d}d_1d_2^2 \si_{d_1}(l)\si_{d_2}\\
&&c_1(L^*)^2=-\hskip-.1in \ma(\sum)_{d_1+d_2=d}\hskip-.1in 
d_2\si_{d_1}(l)\si_{d_2}\quad\mbox{and}\quad
c_1(\wt L^*)^2=-2\hskip-.1in \ma(\sum)_{d_1+d_2=d}\hskip-.1in
d_2\si_{d_1}(l)\si_{d_2}
\eest
\end{lemma}
\pf The relations above follow either by definition, or by aplying
several times (\ref{c1l}) combined with (\ref{x1y}) or (\ref{x2y}) 
(and of course some  simple  algebraic manipulations).  
\medskip

\rem If we distribute  the constraints in Proposition \ref{compp3} 
in all possible ways, formula (\ref{formp3}) becomes: 
\bear\label{3cf}
&\tau_d(p^a,l^b)={2(d-1)(d-2)\over n_jd}\si_d(p^a,l^{b+1})&  \\
\nonumber
&+{2\over n_jd }\ma(\sum)_{d_1=1}^{d-1}\ma(\sum)_{a_1=0}^a
\ma(\sum)_{b_1=0}^b {a\choose a_1}{b\choose b_1}
d_2(2d_1d_2-d)\;\si_{d_1}(p^{a_1},l^{b_1+1})\cdot
\si_{d_2}(p^{a_2},l^{b_2})&
\eear
where in the sum above $d_1+d_2=d$, $a_1+a_2=a$ and $b_1+b_2=b$. 

\step{\bf Example 1} Using a computer program based on (\ref{3cf}) and the 
recursive formulas (\ref{spheres}) for $\si_d$, one recovers for example 
that in $\P^3$ all the degree 2 elliptic invariants are 0 (fact known for 
a very long time) but also one gets new examples, like:

\[ \begin{tabular}{c|ccccc}
& $\tau_3(l^{11})$ && $\tau_5(p,l^{17})$ && $\tau_6(p^{11},l)$ \cr
\hline 
$j\ne 0,1728$ & $6\cdot 25920$ && $6\cdot 15856790593536$ && $6\cdot 13260$
\cr 
$j=0$&  $2\cdot 25920$ && $2\cdot 15856790593536$ && $2 \cdot  13260 $
\cr 
$j=1728$ &  $3\cdot 25920$ && $3\cdot 15856790593536$ && $3\cdot  13260 $
\end{tabular} \]

\step{\bf Example 2} Similarly, when $n=4$, one can use  a computer program
based on the four  steps described in the Section \ref{Comp} to get for
example:
\[ \begin{tabular}{c|ccccc}
& $\tau_3(5H^2,3l,p)$ && $\tau_3(12H^2,l)$   &&$\tau_3(14H^2)$\cr
\hline 
$j\ne 0,1728$ & $6\cdot 42 $ && $6\cdot202680$&& $6\cdot 1305640$
\cr
$j=0$&  $2\cdot 42$ && $2\cdot 202680$ && $2\cdot 1305640 $
\cr
$j=1728$ &  $3\cdot 42$ && $3\cdot 202680$&& $3\cdot 1305640 $
\end{tabular} \]
\medskip

\rem Unfortunately, the number of steps involved in computing the elliptic
invariant $\tau_d$ in $\P^n$ increases extremely fast with $n$. For
example, one can write down the recursive formulas for $n=5,\;6$ that do
not look that complicated (we need to blow up only two strata, $\cz_3$
and $\cz_2$). But the amount of time necessary to run the corresponding 
program is too long to produce interesting examples. 

\section{Appendix}
\renewcommand{\theequation}{A.\arabic{equation}}
\setcounter{equation}{0}
The genus zero perturbed 
invariant and the genus  zero enumerative invariant are equal in $\P^n$
(cf. \cite{rt}), i.e. 
\bear\label{si=phi}
\si_d(H^{j_1},H^{j_2},\dots, H^{j_k})= 
RT_{d,0}(H^{j_1},H^{j_2},H^{j_3}|H^{j_4},\dots, H^{j_k}) 
\eear
 Consequences of Ruan-Tian degeneration formula are:
\bear\label{gendeg}
 RT_{d,1}(\b_1\;|\;\b_2,\dots,\b_l)=
\ma(\sum)_{i_1+i_2=n} \si_d(H^{i_1},H^{i_2},\b_1,\dots,\b_l)
\eear
and that $\si_d$ in $\P^n$ satisfies the following recursive formula:
for  $j_1\ge j_2\ge\dots\ge j_k\ge 2$,
\bear\label{spheres}\nonumber
 \si_d(H^{j_1},H^{j_2},H^{j_3})=
d \si_d(H^{j_1+j_3-1},H^{j_2})-d \si_d(H^{j_1+j_2},H^{j_3-1})
-\si_d(H^{j_1},H^{j_2+1},H^{j_3-1})\\
 \ma(\sum)_{d_1=1}^{d-1}\ma(\sum)_{i=0}^n
\mspa\si_{d_1}(H^{j_1},H^{j_2},H^i)\si_{d_2}(H^{j_3-1},H^{n-i})
\ma(\sum)_{d_1=1}^{d-1}\ma(\sum)_{i=0}^n
\si_{d_1}(H^{j_1},H^{j_3-1},H^i)\si_{d_2}(H^{j_2},H^{n-i})\hskip.1in
\eear
where $\si_d(H^{j_1},H^{j_2},H^{j_3})= 
\si_d(H^{j_1},H^{j_2},H^{j_3},H^{j_4},\dots,H^{j_k})$ and the
conditions $H^{j_4},\dots,H^{j_k}$ are distributed in the right hand side 
in all possible ways.  Note that $\si_1(pt, pt)=1$ gives the initial step 
 of the recursion.

\non Address: Eleny Ionel, MSRI, 1000 Centennial Drive, Berkeley, CA 94720-5070

\non e-mail: {\tt ionel@msri.org}
\end{document}